\documentclass[11pt]{article}

\usepackage{amssymb,latexsym}
\usepackage{amsmath,amscd}
\usepackage{theorem}

\title{Bimodule deformations, Picard groups and contravariant connections}

\author{\textbf{Henrique Bursztyn}\thanks{E-mail: henrique@math.berkeley.edu}
\addtocounter{footnote}{2}\thanks{Research supported by MSRI through NSF grant DMS-9810361 and
    by an Action de Recherche Concertee
    de la Communaut\'e fran\c{c}aise de Belgique.}
  \\[0.1cm]
  D\'epartement de Math\'ematiques,\\
  Universit\'e Libre de Bruxelles,
  CP 218, Campus Plaine\\
  1050 Brussels, Belgium\\[0.5cm]
  \textbf{Stefan Waldmann}\addtocounter{footnote}{2}\thanks{Stefan.Waldmann@physik.uni-freiburg.de}
  \\[0.1cm]
  Fakult{\"a}t f{\"u}r Physik\\
  Albert-Ludwigs-Universit{\"a}t Freiburg\\
  Hermann Herder Stra{\ss}e 3\\
  D 79104 Freiburg\\
  Germany}

\date{July 2002\\[0.5cm] FR-THEP 2002/10}

\newcommand{\im}         {{i}}
\newcommand{\eu}         {{\mathrm e}}

\newcommand{\id}         {{\mathrm {Id}}}
\newcommand{\image}      {{\mathrm {im}}}
\newcommand{\ad}         {{\mathrm {ad}}}
\newcommand{\Aut}        {{\mathrm {Aut}}}
\newcommand{\Inaut}      {{\mathrm {InnAut}}}
\newcommand{\Equiv}      {{\mathrm {Equiv}}}
\newcommand{\Inequiv}    {{\mathrm {InnEquiv}}}
\newcommand{\outequiv}   {{\mathrm {OutEquiv}}}
\newcommand{\Hom}        {{\mathrm {Hom}}}
\newcommand{\End}        {{\mathrm {End}}}
\newcommand{\Der}        {{\mathrm {Der}}}
\newcommand{\PDer}       {{\mathrm {PDer}}}
\newcommand{\IDer}       {{\mathrm {IDer}}}
\newcommand{\st} [1]     {\scriptscriptstyle{{#1}}}
\newcommand{\Ad}         {{\mathrm {Ad}}}

\newcommand{\SP} [1]     {{\left\langle {{#1}} \right\rangle}}

\newcommand{\pr}         {{\mathrm {pr}}}
\newcommand{\cent}       {{\mathcal {Z}}}
\newcommand{\cente}       {{\mathcal {Z}}}
\newcommand{\cl}         {{\mathrm {cl}}}
\newcommand{\curv}         {{\mathrm {curv}}}
\newcommand{\Def}        {{\mathrm {Def}}}
\newcommand{\Defc}        {{\underline{\mathrm {Def}}}}
\newcommand{\Poiss}      {{\mathrm {Poiss}}}
\newcommand{\Bim}        {{\mathrm {EBim}}}
\newcommand{\Bimc}       {{\underline{{\mathrm {EBim}}}}}
\newcommand{\Diff}        {{\mathrm {Diff}}}
\newcommand{\Sympl}        {{\mathrm {Sympl}}}
\newcommand{\Tor}        {{\mathrm {Tor}}}
\newcommand{\GL}        {{\mathrm {GL}}}
\newcommand{\Exp}        {{\mathrm {Exp}}}
\newcommand{\contract}        {{\mathfrak{i}}}
\newcommand{\deriv}        {{\sigma}}

\newcommand{\Conn}        {{\mathrm {Conn}}}
\newcommand{\Connc}        {{\underline{\mathrm {Conn}}}}

\newcommand{\A}          {\mathcal{A}}
\newcommand{\B}          {\mathcal{B}}
\newcommand{\qA}         {\boldsymbol{\mathcal{A}}}
\newcommand{\qB}         {\boldsymbol{\mathcal{B}}}

\newcommand{\Pic}         {\mathrm{Pic}}
\newcommand{\Picc}        {\underline{\mathrm{Pic}}}
\newcommand{\PicZ}         {\mathrm{Pic}_{\st{\cent}}}
\newcommand{\PicA}         {\mathrm{Pic}_{\st{\A}}}

\newcommand{\PiccZ}        {\underline{\mathrm{Pic}}_{\st{\mathcal{Z}}}}
\newcommand{\PiccA}        {\underline{\mathrm{Pic}}_{\st{\A}}}

\newcommand{\LogInnAut}  {\mathrm{LogInnAut}}
\newcommand{\InnDer}     {\mathrm{InnDer}}
\newcommand{\HdR}        {H_{\mathrm{\scriptscriptstyle dR}}}
\newcommand{\LocInnDer}  {\mathrm{LocInnDer}}

\newcommand{\M}          {\mathrm{M}}
\newcommand{\qM}         {\boldsymbol{\M}}

\newcommand{\qP}         {\boldsymbol{P}}

\newcommand{\X}          {\mathrm{X}}
\newcommand{\Y}          {\mathrm{Y}}

\newcommand{\AXA}        {_{\scriptscriptstyle{\A}}{\X}_{\scriptscriptstyle{\A}}}
\newcommand{\BXA}        {_{\scriptscriptstyle{\B}}{\X}_{\scriptscriptstyle{\A}}}

\newcommand{\qX}         {\boldsymbol{\X}}
\newcommand{\qY}         {\boldsymbol{\Y}}

\newcommand{\tensorA}    {\otimes_{\scriptscriptstyle{\A}}}
\newcommand{\tensorB}    {\otimes_{\scriptscriptstyle{\B}}}

\newcommand{\modA}       {_{\scriptscriptstyle{\A}}\mathfrak{M}} 
\newcommand{\modB}       {_{\scriptscriptstyle{\B}}\mathfrak{M}} 

\newcommand{\funcX}      {\mathfrak{F}_{\scriptscriptstyle{\X}}}
\newcommand{\func}       {\mathfrak{F}}

\newtheorem{lemma} {Lemma} [section]
\newtheorem{proposition} [lemma] {Proposition}

\newtheorem{theorem} [lemma] {Theorem}
\newtheorem{corollary} [lemma] {Corollary}
\newtheorem{definition}[lemma] {Definition}

\theorembodyfont{\rm} 

\newtheorem{example}[lemma] {Example}

\newtheorem{remark}[lemma]{Remark}

\newenvironment{proof}{{\sc Proof:}}{{\hspace*{\fill} $\square$\\}}

\numberwithin{equation}{section}

\begin{document}

\maketitle

\begin{abstract}
    We study deformations of invertible bimodules and the behavior of
    Picard groups under deformation quantization.  While 
    $K_0$-groups are known to be stable under formal deformations of
    algebras, Picard groups may change drastically.  We identify the
    semiclassical limit of bimodule deformations as contravariant
    connections and study the associated deformation quantization
    problem.  Our main focus is on formal deformation quantization of
    Poisson manifolds by star products.
\end{abstract}


\section{Introduction}

Some important mathematical formulations of quantization are based on
replacing a classical algebra of observables $\A$ by a noncommutative
one, obtained from $\A$ through a deformation procedure \cite{BFFLS78}
(see also \cite{Land98,Rief89} for more analytical approaches).  It is
natural to study the behavior of algebraic invariants, such as
$K$-theory and Hochschild (co)homology, in this process, see e.g.
\cite{NT95a,NT95b,Ros96,WX98}.  Building on \cite{BuWa2002}, the main
objective of this paper is to discuss the behavior of Picard groups
under formal deformation quantization.

More precisely, let $\A$ be a commutative unital algebra over a
commutative, unital ring $k$, and let $\qA=(\A[[\lambda]],\star)$ be
a formal deformation of $\A$ \cite{Ger64}.  The classical limit map
\begin{equation}
    \label{eq:clintro}
    \cl : \qA \longrightarrow \A, 
    \quad 
    \sum_{r=0}^\infty \lambda^r a_r \mapsto a_0,
\end{equation}
induces a group homomorphism
\begin{equation}
    \label{eq:picintro}
    \cl_*: \Pic(\qA) \longrightarrow \Pic(\A),
\end{equation}
and our goal is to study the kernel and image of $\cl_*$, with special
focus on the case where the deformations come from star products on
Poisson manifolds.  Recall that the Picard group of an arbitrary
unital algebra is defined as the group of isomorphism classes of
invertible bimodules over it (see Definition~\ref{def:ebimod}); for example,
if $\A = C^\infty(M)$ is the algebra of complex-valued smooth
functions on a manifold $M$, then $\Pic(\A)$ is the semi-direct
product
$$
\Pic(C^\infty(M)) = \Diff(M)\ltimes \Pic(M),
$$
where $\Diff(M)$ is the diffeomorphism group and $\Pic(M)\cong
H^2(M,\mathbb{Z})$ is the \emph{geometric} Picard group, consisting of
isomorphism classes of complex line bundles over $M$; the action is
given by pull-back of line bundles.

The behavior of $K_0$-groups under deformations was discussed in
\cite{Ros96,BuWa2000b}: the classical limit map (\ref{eq:clintro})
induces a map of $K_0$-groups,
$$
K_0(\qA) \longrightarrow K_0(\A),
$$
which turns out to be a group isomorphism \cite[Thm.~4]{Ros96}.
In the case of (\ref{eq:picintro}), we will see that $\cl_*$ 
is far from being an isomorphism in general.

In the purely algebraic setting, the kernel of $\cl_*$ is in
one-to-one correspondence with outer self-equivalences of $\qA$, i.e.,
the group of automorphisms of $\qA$ of the form $T = \id + O(\lambda)$
modulo inner automorphisms. The image of $\cl_*$ is described in terms
of a canonical action $\Phi$ of $\Pic(\A)$ on the moduli space of
equivalence classes of formal deformations of $\A$ \cite{Bu2002a}.
Roughly speaking, two deformations $\star'$ and $\star$ of $\A$ are in
the same $\Phi$-orbit if and only if there exists an
$(\A,\A)$-invertible bimodule $\X$ that can be deformed into a
$(\star',\star)$-equivalence bimodule. In fact, the $\Phi$-orbits
characterize Morita equivalent deformations of $\A$.

Just as associative algebraic deformations of $\A$ correspond to
Poisson brackets, the semiclassical limit of a
$(\star',\star)$-bimodule deformation of an 
$(\A,\A)$-invertible bimodule $\X$ is an algebraic generalization of what is
known in Poisson geometry as a linear \emph{contravariant connection}
\cite{Rui2000,Hueb90,Vais91}. We note that contravariant connections on $\X$
arising in the semiclassical limit of a $(\star',\star)$-bimodule
deformation must have a fixed curvature $\tau$, depending only on
$\star'$ and $\star$. More generally, we discuss the existence and
classification of bimodule deformations ``in the direction'' of a
fixed contravariant connection, in the same spirit of deformation
quantization of Poisson algebras.

In the geometric context of star products on Poisson manifolds, our
results are as follows. We consider star products satisfying two extra
conditions (see C0) and C1) in Sections \ref{sec:innerauto} and
\ref{sec:linedef}) relating the Poisson center and the Poisson
derivations to the center and derivations of the star product.  These
conditions are always satisfied, e.g.,\ by symplectic star products
and, on arbitrary Poisson manifolds, by the star products constructed
in \cite{CFT01,CFT02}.  For such star products on a Poisson manifold
$(M,\pi)$, we show that, for a fixed line bundle with contravariant
connection satisfying the curvature condition, there always exist
corresponding bimodule deformations. Moreover, equivalence classes of
bimodule deformations in the direction of a given contravariant
connection are in bijection with $H^1_{\pi}(M,\mathbb{C})[[\lambda]]$,
the set of formal power series with coefficients in the first Poisson
cohomology group.

Under suitable conditions, always satisfied by symplectic star products,
we show that the kernel of $\cl_*$  is in bijection with
\begin{equation}
    \label{eq:kerintro}
    \frac{H^1_{\pi}(M,\mathbb{C})}
    {2\pi \im H^1_{\pi}(M,\mathbb{Z})}
    + \lambda H^1_{\pi}(M,\mathbb{C})[[\lambda]].
\end{equation}
In particular, $\cl_*$ is injective if and only if
$H^1_{\pi}(M,\mathbb{C})=0$.

Finally, we give a description of the image of $\cl_*$ for symplectic
star products.  As an example, we show that if
$[\psi^*(\star)]=[\star]$ for all symplectomorphisms $\psi$ (here
$[\star]$ denotes the equivalence class of the star product $\star$),
then we can write the image of $\cl_*$ as the semi-direct product
$$
\cl_*(\Pic(\qA)) = \Sympl(M)\ltimes \Tor(H^2(M,\mathbb{Z})),
$$
where $\Sympl(M)$ denotes the group of symplectomorphisms of $M$
and\ $\Tor(H^2(M,\mathbb{Z}))$ is the subgroup of torsion elements in
$H^2(M,\mathbb{Z})$; the action is given by pull-back.  This is the
case for arbitrary star products on Riemann surfaces and on
$\mathbb{C}P^n$.  We also construct explicit examples where the image
of the map $\cl_*^r = \pr \circ \cl_*$,
$$
\Pic(\qA) \stackrel{\cl_*}{\longrightarrow} 
\Pic(\A)=\Diff(M)\ltimes H^2(M,\mathbb{Z}) 
\stackrel{\pr}{\longrightarrow} H^2(M,\mathbb{Z}) 
$$
contains non-torsion elements, but it seems hard to describe
exactly how big this image can be. Nevertheless, we show that
$\cl_*^r$ is onto if and only if $H^2(M,\mathbb Z)$ only contains
torsion elements.

The paper is organized as follows.  In Section~\ref{sec:prelim}, we
review the basics of equivalence bimodules and Picard groups.  In
Section~\ref{sec:deform}, we define the map $\cl_*$ and give
descriptions of its kernel and image in purely algebraic terms.
Section~\ref{sec:bidef} discusses bimodule deformations and the
contravariant connections in their semiclassical limit.
Section~\ref{sec:auto} collects results on the structure of the
automorphism group of star product algebras, used in
Section~\ref{sec:linedef} to discuss deformation quantization of line
bundles over Poisson manifolds. Finally, Section~\ref{sec:kerimagecl}
describes the kernel and the image of $\cl_*$ in cohomological terms
for the case of star products.

\noindent{\bf Acknowledgements:}
    We would like to thank M.~Bordemann,
    M.~Crainic, A.~Gammella, N.~Giulini, S.~Gutt, R.~Fernandes, D.~Martinez, A.~Serra
    and A.~Weinstein for valuable discussions and helpful advice. H.B.
    thanks MSRI and ULB for their hospitality while this work was being
    developed.


\section{Preliminaries on equivalence bimodules  and Picard groups}
\label{sec:prelim}

In this section we recall definitions and standard facts about
equivalence bimodules and Picard groups. The reader should consult
\cite{Bass68} for details. 

\subsection{Equivalence bimodules and Picard groups}
\label{sec:equipic}

Let $k$ be a commutative unital ring of characteristic zero, and let
$\A$ and $\B$ be associative unital algebras over $k$.  For
convenience we shall assume $\mathbb{Q} \subseteq k$.
\begin{definition}
    \label{def:ebimod}
    A $(\B,\A)$-bimodule $\X$ is called an \emph{equivalence bimodule}
    if there exists an $(\A,\B)$-bimodule $\X^*$ so that 
    $\X\tensorA\X^* \cong \B$ as $(\B,\B)$-bimodules and 
    $\X^* \tensorB \X \cong \A$ as $(\A,\A)$-bimodules.
\end{definition}
Let $\Bimc(\B,\A)$ denote the category of $(\B,\A)$-equivalence
bimodules, with bimodule homomorphisms as morphisms. The set of
isomorphism classes of objects in $\Bimc(\B,\A)$ is denoted by
$\Bim(\B,\A)$. When $\A = \B$, we denote $\Bimc(\A,\A)$ by
$\Picc(\A)$, and objects in $\Picc(\A)$ are also called
$(\A,\A)$-\emph{invertible bimodules}. In this case, the set
$\Bim(\A,\A)$ forms a group under $\tensorA$, denoted by $\Pic(\A)$.
\begin{definition}
    \label{def:picgr}
    The group $\Pic(\A)$ is called the \emph{Picard group} of $\A$.
\end{definition}
We denote the isomorphism class of an object $\X$ in $\Bimc(\B,\A)$ by
$[\X] \in \Bim(\B,\A)$.

We recall a few facts from Morita theory \cite{Morita58}, see also
e.g.~\cite{Lam99}.
\begin{enumerate}
\item Let $\modA$ (resp.\  $\modB$) denote the category of left
    $\A$-modules (resp.\  left $\B$-modules). Each object $\X$ in
    $\Bimc(\B,\A)$ defines a functor
    \begin{equation}
        \label{eq:functor} 
        \funcX = \X\tensorA \cdot : {\modA} \longrightarrow {\modB},
        \quad
        \funcX(\M) = \X \tensorA \M,
    \end{equation}
    which is an equivalence of categories. Conversely, any equivalence
    functor $\func: {\modA} \longrightarrow {\modB}$ is isomorphic to
    a functor of the form (\ref{eq:functor}), for some 
    $\X \in \Bimc(\B,\A)$. When $\A = \B$, this correspondence induces
    an isomorphism between $\Pic(\A)$ and the group of isomorphism
    classes of self-equivalence functors of $\modA$, with group
    operation given by composition.
    
\item If $\X \in \Bimc(\B,\A)$, then
    $[\X^*]=[\Hom_{\st{\A}}(\X_{\scriptscriptstyle{\A}},\A)]$ and we
    have an algebra isomorphism $\End_{\st{\A}}(\X_{\st{\A}})\cong\B$
    via the left $\B$-action.  In particular, this action is
    non-degenerate: $bx =0$ for all $x \in \X$ implies that $b=0$.
    Similar results hold for the right action of $\A$.
    
\item Any $\X \in \Bimc(\B,\A)$ is finitely generated, projective and
    full both as a left $\B$-module and right $\A$-module. In
    particular, as a right $\A$-module, $\X_{\st{\A}}\cong P\A^n$,
    where $P\in M_n(\A)$ is a full idempotent and $n \in \mathbb{N}$;
    fullness means that $M_n(\A) P M_n(\A) = M_n(\A)$. As a result of
    Morita's theorem \cite{Morita58}, $\X \in \Bimc(\B,\A)$ if and
    only if there exists a full idempotent $P\in M_n(\A)$ such that
    $\X_{\st{\A}}\cong P\A^n$ as right $\A$-modules and $PM_n(\A)P
    \cong \B$. This shows in particular that $\Bim(\B,\A)$ and
    $\Pic(\A)$ are indeed sets.
\end{enumerate}

\subsection{The contribution of algebra automorphisms}
\label{sec:contribauto}

We denote the group of automorphisms of $\A$ by $\Aut(\A)$, and
$\Inaut(\A)$ denotes the normal subgroup of inner automorphism, i.e.,
automorphisms of the form $\Ad(u): a \mapsto u a u^{-1}$ for an
invertible element $u \in \A$.

Let $\X \in \Bimc(\B,\A)$, $\phi \in \Aut(\B)$ and $\psi \in
\Aut(\A)$. We define a new element 
$$
_{\st{\phi}}\X_{\st{\psi}} \in
\Bimc(\B,\A)
$$ 
by setting $_{\st{\phi}}\X_{\st{\psi}} \cong \X$ as an
additive group and changing the left and right actions by
$$
b \cdot x := \phi(b)x \qquad \mbox{ and } \qquad x \cdot a := x\psi(a),
$$
for $b \in \B$, $a \in \A$ and $x \in \X$. In this notation, 
$\X = {_{\st{\id}}}\X_{\st{\id}}$.  It is easy to check that
\begin{equation}
    \label{eq:Atensor}
    _{\st{\phi}}\X_{\st{\psi}} 
    \cong {_{\st{\phi}}}\B \tensorB \X \tensorA \A_{\st{\psi}}.
\end{equation}
\begin{proposition}
    \label{prop:transit}
    Let $\X, \Y \in \Bimc(\B,\A)$ be such that $\X_{\st{\A}} \cong
    \Y_{\st{\A}}$ (as right $\A$-modules).  Then there exists $\phi
    \in \Aut(\B)$ so that ${_{\st{\phi}}}\X_{\st{\id}} \cong \Y$ (as
    $(\B,\A)$-bimodules).  Moreover, ${_{\st{\phi}}}\X_{\st{\id}}
    \cong \X$ as $(\B,\A)$-bimodules if and only if $\phi \in
    \Inaut(\B)$.
\end{proposition}
\begin{proof}
    Let $j:\X_{\st{\A}}\longrightarrow \Y_{\st{\A}}$ be a right
    $\A$-module isomorphism. Fix $b \in \B$ and consider the map from
    $\X$ to itself defined by
    $$
    x \quad\mapsto\quad j^{-1}(bj(x)).
    $$
    Since this map is right $\A$-linear and $\B \cong
    \End_{\st{\A}}(\X_{\st{\A}})$, it defines a unique element
    $\tilde{b} \in \B$ so that $\tilde{b}\cdot x =j^{-1}(bj(x))$.
    This procedure defines a map
    $$
    \phi:\B \longrightarrow \B, \;\; \phi(b)=\tilde{b},
    $$
    which is an automorphism of $\B$.  From
    the equation $j(\phi(b)x)=bj(x)$, it follows that
    $$
    j:{_{\st{\phi}}\X_{\st{\id}}}\longrightarrow \Y
    $$
    is a bimodule isomorphism.

    For the second statement first assume $\phi = \Ad(u)$ to be inner.
    Then $x \mapsto u^{-1}x$ is the desired bimodule isomorphism.
    Conversely, if $f: \X \to {_{\st{\phi}}}\X_{\st{\id}}$ is a
    bimodule isomorphism, then $f(x) = ux$ for some $u \in \B$ since
    $\B \cong \End_{\st{\A}}(\X_{\st{\A}})$. As the action of $\B$ is
    non-degenerate and $f$ is invertible, $u$ must also be invertible.
    Finally,
    $$
    ubx = f(bx) = \phi(b) f(x) = \phi(b) u x \; \implies \; 
    \phi(b) = ubu^{-1}. 
    $$
\end{proof}

In order to see how much $\Aut(\A)$ contributes to $\Pic(\A)$, we
consider the map
\begin{equation}
    \label{eq:ldef}
    l: \Aut(\A) \longrightarrow \Pic(\A), 
    \quad l(\phi)=[{_{\st{\phi^{-1}}}}\A_{\st{\id}}].
\end{equation}
A simple computation shows that, for $\phi, \psi \in \Aut(\A)$,
${_{\st{\phi^{-1}}}\A_{\st{\id}}} \tensorA
{_{\st{\psi^{-1}}}\A_{\st{\id}}} \cong
{_{\st{\psi^{-1}\phi^{-1}}}}\A_{\st{\id}}$.  
So the map $l$ is a group homomorphism.  By
Proposition~\ref{prop:transit}, we have an exact sequence of groups
\begin{equation}
    \label{eq:exactseq}
    1 \to \Inaut(\A) \to \Aut(\A) \stackrel{l}{\to} \Pic(\A).
\end{equation}
Thus the image of $l$ is a subgroup of $\Pic(\A)$ isomorphic to
$\Aut(\A)/\Inaut(\A)$, the group of \emph{outer automorphisms} of $\A$.

\subsection{Remarks on the center} 
\label{sec:center}

Let $\cente(\A)$ denote the center of $\A$. If $\X \in \Picc(\A)$, it
is not generally true that $\cente(\A)$ acts the same way on the left
and right of $\X$. Let $\PiccZ(\A)$ be the category of
$(\A,\A)$-invertible bimodules satisfying
$$
z x = x z,
$$
for all $x \in \X$ and $z \in \cente(\A)$ and let $\PicZ(\A)$ be the
corresponding group of isomorphism classes of objects.
\begin{proposition}\label{prop:center}
    There exists a group homomorphism 
    $h:\Pic(\A) \longrightarrow \Aut(\cente(\A))$
    so that
    \begin{equation}
        \label{eq:hexactsequence}
        1 \to \PicZ(\A) \to \Pic(\A) \stackrel{h}{\to} \Aut(\cente(\A))
    \end{equation}
    is exact. Moreover, if $\A$ is commutative then $h$ is split by
    the map $l : \Aut(\A) \to \Pic(\A)$ as in \eqref{eq:ldef}.
\end{proposition}
\begin{proof}
    Let $\X \in \Picc(\A)$.
    If $z \in \cente(\A)$, then the map $x \mapsto xz$ is a bimodule
    endomorphism of $\X$. Since $\End_{\st{\A}}(\X_{\st{\A}}) \cong \A$,
    there exists a unique $a=a(z) \in \A$ with $ax = xz$. It is clear that
    $a \in \cente(\A)$ and that we have an induced algebra homomorphism
    $$
    \phi_{\X}:\cente(\A) \longrightarrow \cente(\A), \;\; z \mapsto a(z).
    $$
    If $\X, \Y \in \Picc(\A)$ and $x \in \X, y\in \Y$, then 
    $$
    x\otimes y z = x \otimes \phi_{\Y}(z) y = x \phi_{\Y}(z) \otimes y = 
    \phi_{\X}(\phi_{\Y}(z))x \otimes y.
    $$
    Hence $\phi_{\X\tensorA \Y} = \phi_{\X} \phi_{\Y}$. 
    Since $\phi_{\A} = \id$, it follows that
    $\phi_{\X} \phi_{\X^*}= \phi_{\X^*} \phi_{\X}  =\id$ 
    and $\phi_{\X}$ is an automorphism of $\cente(\A)$.
    Thus 
    $$
    h:\Pic(\A) \longrightarrow \Aut(\cente(\A)),\;\;
    \X \mapsto \phi_{\X}
    $$ 
    is a homomorphism.
    The result now easily follows.
\end{proof}

Let us now assume that $\A$ is \emph{commutative}. In this case, we
regard $\Aut(\A)$ as a subgroup of $\Pic(\A)$ via the embedding 
$l: \Aut(\A) \to \Pic(\A)$, $\phi \mapsto [{_{\st{\phi^{-1}}}}\A]$.
\begin{corollary}\label{cor:semidir}
    If $\A$ is commutative, then $\PicA(\A)$ is a normal subgroup of 
    $\Pic(\A)$, and $\Pic(\A)$ is the semi-direct product of
    $\PicA(\A)$ and $\Aut(\A)$. The induced product on 
    $\Aut(\A) \times \PicA(\A)$ is given by 
    \begin{equation} 
        \label{eq:operation}
        (\phi,[\X])\cdot (\psi,[\Y]) = 
        (\phi\circ \psi, [{_{\st{\psi}}}\X_{\st{\psi}} \tensorA \Y]).
    \end{equation}
\end{corollary}
\begin{proof}
    It follows from Proposition~\ref{prop:center} that $\PicA(\A)$ is
    a normal subgroup of $\Pic(\A)$, and
    Proposition~\ref{prop:transit} implies that the map
    $$
    f : \Aut(\A)\times \PicA(\A) \longrightarrow \Pic(\A), \;\;
    f(\phi,[\X]) = [{_{\st{\phi^{-1}}}}\X],
    $$ 
    is a bijection. Finally, note that 
    $$
    {_{\st{\phi^{-1}}}\X} \tensorA {_{\st{\psi^{-1}}}\Y} \cong 
    {_{\st{\phi^{-1}}}\A} \tensorA \X \tensorA{_{\st{\psi^{-1}}}\A} \tensorA \Y \cong
    {_{\st{\phi^{-1}}}\A} \tensorA \X \tensorA {\A_{\st{\psi}}} \tensorA \Y.
    $$
    On the other hand,
    $$
    \X \tensorA {\A_{\st{\psi}}} 
    \cong
    \X_{\st{\psi}}
    \cong 
    {_{\st{\psi^{-1}}}}\A \tensorA {_{\st{\psi}}\X_{\st{\psi}}}.
    $$
    Therefore
    $$
    {_{\st{\phi^{-1}}}\X} \tensorA {_{\st{\psi^{-1}}}\Y} \cong
    {_{\st{\psi^{-1}}\st{\phi^{-1}}}}\A \tensorA {_{\st{\psi}}\X_{\st{\psi}}} \tensorA \Y,
    $$
    proving (\ref{eq:operation}).
\end{proof}

For a commutative algebra $\A$, the group $\PicA(\A)$ is called the
\emph{classical} or \emph{commutative} Picard group of $\A$.
\begin{example}
    \label{ex:semidir}
    Let $M$ be a smooth manifold and let $\A = C^\infty(M)$ be the
    algebra of smooth complex-valued functions on $M$. Then
    \begin{equation}
        \label{eq:PicM}
        \PicA(\A) \cong \Pic(M) \cong H^2(M,\mathbb{Z}).
    \end{equation}
    The first identification is a result of Serre-Swan's theorem
    (see e.g.~\cite{Bass68}), where we identify line bundles with
    invertible bimodules by $L \mapsto \Gamma^\infty(L)$;
    the second is the Chern class map
    (see e.g.~\cite{Hirz95}). Hence  $\Pic(\A)$ is the semi-direct product of
    $H^2(M,\mathbb{Z})$ and $\Diff(M)$, where $\Diff(M)$ is the group
    of diffeomorphisms of $M$ acting on $H^2(M,\mathbb{Z})$ by
    pull-back:
    $$
    (\phi_1,l_1)\cdot (\phi_2,l_2) 
    = (\phi_1\circ \phi_2, \phi_2^*(l_1) + l_2),
    $$
    where $\phi_1,\phi_2 \in \Diff(M)$, $l_1,l_2 \in H^2(M,\mathbb{Z})$.
\end{example}


\section{Picard groups of deformed algebras}
\label{sec:deform}

\subsection{Deformations of algebras and modules}
\label{sec:algmoddef}

We will recall the basic definitions of formal deformations of
associative algebras and modules over them, see \cite{GS88}.

Let $\A$ be a commutative and associative unital algebra over $k$,
with $\mathbb{Q} \subseteq k$. We note that many of the definitions
and results to follow hold for noncommutative algebras as well, but we
will restrict ourselves to the commutative setting. Let
$\A[[\lambda]]$ denote the space of formal power series with
coefficients in $\A$ and formal parameter $\lambda$.
\begin{definition}
    \label{def:algdef}
    A formal \emph{deformation} of $\A$ is a $k[[\lambda]]$-bilinear
    associative product $\star$ on $\A[[\lambda]]$ of the form
    \begin{equation}\label{eq:starp}
        a_1 \star a_2 = \sum_{r=0}^{\infty} \lambda^r C_r(a_1,a_2),
        \;\; a_1, a_1 \in \A, 
    \end{equation}
    where $C_r :\A \times \A \longrightarrow \A$ are bilinear maps,
    and $C_0(a_1,a_2) = a_1 a_2$. We also require that $1 \star a = a
    \star 1 = a$, for all $a \in \A$, where $1 \in \A$ is the unit
    element.
\end{definition}
We denote the resulting deformed $k[[\lambda]]$-algebra by
$\qA = (\A[[\lambda]],\star)$. 
\begin{definition}
    \label{def:algequiv}
    Two deformations $\qA_1=(\A[[\lambda]],\star_1) $ and 
    $\qA_2=(\A[[\lambda]],\star_2)$
    are \emph{equivalent} if there are $k$-linear maps 
    $T_r:\A \longrightarrow \A$ so that 
    $$
    T = \id + \sum_{r=1}^\infty \lambda^r T_r: 
    \; \qA_1 \longrightarrow \qA_2
    $$
    is an algebra isomorphism. 
\end{definition}

The equivalence class of $\star$ is denoted by $[\star]$. The set of
equivalence classes of deformations of $\A$ is denoted by $\Def(\A)$.
A deformation $\star=\sum_{r=0}^\infty \lambda^r C_r$ defines a
Poisson bracket on $\A$ by
\begin{equation}
    \label{eq:PoissonBracketDef}
    \{a_1,a_2\} := C_1(a_1,a_2) - C_1(a_2,a_1),
\end{equation}
and equivalent deformations define equal Poisson brackets (see
e.g.~\cite{SilWein99}).  We denote the set of equivalence classes of
deformations corresponding to a fixed Poisson bracket on $\A$ by
$\Def(\A,\{\,,\,\})$.

We note that $\Aut(\A)$ acts on formal deformations of $\A$ by
$\star \mapsto \hat{\star} = \psi^*(\star)$, where 
$$
a_1 \mathbin{\hat{\star}} a_2 = \psi^{-1}(\psi(a_1)\star \psi(a_2)), 
\quad \psi \in \Aut(\A).
$$
It is simple to check that two deformations $\qA_1$ and $\qA_2$ are
isomorphic as $k[[\lambda]]$-algebras if and only if there exists a
$\psi \in \Aut(\A)$ with $[\psi^*(\star_1)] = [\star_2]$.  Note,
however, that this action of $\Aut(\A)$ typically changes the Poisson
bracket as defined in (\ref{eq:PoissonBracketDef}).

The following definition of a star product \cite{BFFLS78}  provides
our main example of formal deformations.
\begin{definition}
    A \emph{star product} on a manifold $M$ is a formal deformation
    $\star = \sum_{r=0}^\infty\lambda^r C_r$ of $C^\infty(M)$ in the
    sense of Definition~\ref{def:algdef} so that each $C_r$ is a
    bidifferential operator.
\end{definition}
In this case, $\Def(M)$ denotes the moduli space of equivalence
classes of star products on $M$; the subspace of equivalence classes
of star products corresponding to a fixed Poisson structure $\pi$ on
$M$ is denoted by $\Def(M,\pi)$. Here $\pi \in
\Gamma^\infty(\bigwedge^2TM)$ is the Poisson bivector field such that
$\{f,g\} = \pi(df, dg)$ is the Poisson bracket, see
e.g.~\cite{Vais94,SilWein99} for details on Poisson geometry.  For
star products, we adopt the convention
\begin{equation}
    \label{eq:convention}
    \{f,g\} := \frac{1}{\im} (C_1(f,g) - C_1(g,f)), 
    \quad f,g \in C^\infty(M).
\end{equation}

We shall now consider module deformations of a
right module $\M_{\st{\A}}$ over an arbitrary commutative unital 
algebra $\A$.
\begin{definition}
    \label{definition:defmodule}
    A \emph{deformation} of $\M_{\st{\A}}$ with respect to a
    deformation $\star$ of $\A$ (or, in short, a $\star$-module
    deformation of $\M_{\st{\A}}$) is a $k[[\lambda]]$-bilinear right
    $\qA$-module structure $\bullet$ on $\M[[\lambda]]$ of the form
    \begin{equation}
        m \bullet a = 
        \sum_{r=0}^\infty \lambda^r R_r(m,a), 
        \;\; m \in \M, \; a \in \A, 
    \end{equation}
    where each $R_r: \M \times \A \longrightarrow \M$ is bilinear and
    $R_0(m,a)=ma$ is the original module structure.
\end{definition}
We denote the deformed module $(\M[[\lambda]],\bullet)$ by
$\qM_{\st{\qA}}$. Deformations of left modules are
defined analogously.
\begin{definition}
    \label{def:modequiv}
    Two $\star$-module deformations $\qM_1$ and $\qM_2$ of $\M$
    are called \emph{equivalent} if
    there exist $k$-linear maps $T_r: \M \longrightarrow \M$ so that
    \begin{equation}
        \label{eq:equivDefMod}
        T = \id + \sum_{r=1}^\infty\lambda^r T_r: \;
        \qM_1 \longrightarrow \qM_2
    \end{equation}
    is an $\qA$-module isomorphism.
\end{definition}

As in the case of deformations of algebras, the automorphisms
$\Aut(\M_{\st{\A}})$ act on $\star$-module deformations of
$\M_{\st{\A}}$ by $\hat{\bullet}= \psi^*(\bullet)$, where
\begin{equation}
    \label{eq:hatbullet}
    m \mathbin{\hat{\bullet}} a = \psi^{-1}(\psi(m) \bullet a), 
    \quad m \in \M, \; a \in \A.
\end{equation}
Two module deformations $\qM_1$ and $\qM_2$ of $\M_{\st{\A}}$ with
respect to $\star$ are isomorphic if and only if there exists $\psi
\in \Aut(\M_{\st{\A}})$ so that $\bullet_1$ is equivalent to
$\psi^*(\bullet_2)$.

Recall from \cite{Ros96,BuWa2000b}
that if $\M_{\st{\A}}$ is finitely generated and projective, then
module deformations always exist with respect to any
deformation $\star$ and are still finitely generated and projective.
Moreover, they are unique, up to equivalence.
On the other hand, any finitely generated projective module over a
deformed algebra $\qA=(\A[[\lambda]],\star)$ is isomorphic to a
deformation of an $\A$-module with respect to $\star$.
Hence the classical limit map 
$$
\cl: \qA \longrightarrow \A, \;\; \sum_{r=0}^\infty a_r 
\lambda_r \mapsto a_0
$$
induces a map of $K_0$-groups
\begin{equation}
    \label{eq:K0}
    \cl_* : K_0(\qA) \longrightarrow K_0(\A), 
    \quad [\qM_{\st{\qA}}] \mapsto [\M_{\st{\A}}],
\end{equation}
where $\qM_{\st{\qA}} \cong (\M_{\st{\A}}[[\lambda]],\bullet)$,
and this map is a group isomorphism \cite[Thm.~4]{Ros96}.

\subsection{Bimodule deformations and the classical-limit map}
\label{sec:bidefcl}

Let $\BXA$ be a $(\B,\A)$-bimodule, and consider deformations
$\qB = (\B[[\lambda]],\star')$ and $\qA=(\A[[\lambda]],\star)$.
\begin{definition}
    \label{def:bimoddef}
    A \emph{bimodule deformation} of $\BXA$ with respect to $\star'$
    and $\star$ (or, in short, a $(\star',\star)$-bimodule deformation
    of $\BXA$) is a $(\qB,\qA)$-bimodule structure on $\X[[\lambda]]$
    (with right and left module structures being
    $k[[\lambda]]$-bilinear) of the form
    \begin{equation}
        \label{eq:bimoddef}
        b \bullet' x 
        = \sum_{r=0}^\infty \lambda^r R_r'(b,x),
        \quad
        x \bullet a 
        = \sum_{r=0}^\infty \lambda^r R_r(x,a), 
    \end{equation}
    for $a \in \A$, $b \in \B$ and $x \in \X$.  Here $R_r':\B \times
    \X \longrightarrow \X$ and $R_r: \X \times \A \longrightarrow \X$
    are bilinear and so that $R'_0(b,x) = b x$ and $R_0(x,a) = x a$
    are the original module structures.
\end{definition}
We denote a $(\star',\star)$-bimodule deformation of $\BXA$ by
$(\X[[\lambda]],\bullet',\bullet) = {_{\st{\qB}}}\qX_{\st{\qA}}$.  For
a fixed deformation $\star$, we will later discuss conditions on
$\star'$ guaranteeing the corresponding bimodule deformations to
exist.

Analogously to the case of algebras and modules, bimodule
automorphisms $\psi$ of $\BXA$ act on $(\star',\star)$-bimodule
deformations by $(\X[[\lambda]],\bullet',\bullet) \mapsto
(\X[[\lambda]],\hat{\bullet}',\hat{\bullet})$, where
$$
b \mathbin{\hat{\bullet}}' x := \psi^{-1}(b \bullet' \psi(x))
\quad\textrm{and}\quad
x \mathbin{\hat{\bullet}} a := \psi^{-1} (\psi(x)\bullet a).
$$
Two $(\star',\star)$-bimodule deformations of $\BXA$ are isomorphic
if and only if they are equivalent up to this action.

Let us now assume that
$\qA = (\A[[\lambda]],\star)$ and $\qB = (\B[[\lambda]],\star')$ are
Morita equivalent. 
\begin{proposition}
    \label{prop:bimdef}
    Any object ${\qX} \in \Bimc(\qB,\qA)$ is isomorphic to a
    $(\star',\star)$-bimodule deformation of some element
    $\X\in\Bimc(\B,\A)$, and this gives rise to a well-defined map
    \begin{equation}
        \label{eq:clpic}
        \cl_*: \Bim(\qB,\qA) 
        \longrightarrow \Bim(\B,\A),\quad [\qX ] \mapsto [\X].
    \end{equation}
\end{proposition}
\begin{proof}
    Let ${_{\st{\qB}}}\qX_{\st{\qA}} \in \Bimc(\qB,\qA)$. As a right
    $\qA$-module, $\qX_{\st{\qA}}$ is isomorphic to $\qP\star \qA^n$,
    for some full idempotent $\qP = P + O(\lambda) \in M_n(\qA)$ with
    an idempotent $P \in M_n(\A)$. In particular, $\qB \cong \qP \star
    M_n(\qA) \star \qP$.

    Since $\qP \star \qA^n \cong P\A^n[[\lambda]]$ as
    $k[[\lambda]]$-modules, see e.g.~\cite{BuWa2000b}, we have an induced
    $(\qB,\qA)$-bimodule structure on $P\A^n[[\lambda]]$, which is a
    bimodule deformation of $\Y = P\A^n \in \Bimc(\B,\A)$ with respect
    to $\star'$ and $\star$ (recall that the fullness of $\qP$ implies
    that $P$ is automatically full). We denote this bimodule
    deformation by $\qY$.

    Since $\qX$ is isomorphic to $\qY$ as a right $\qA$-module, it
    follows from Proposition~\ref{prop:transit} that there exists
    $\psi = \psi_0 + O(\lambda) \in \Aut(\qB)$ such that
    ${_{\st{\qB}}}\qX_{\st{\qA}}$ is isomorphic to
    ${_{\st{\psi}}}\qY_{\st{\qA}}$, which is a deformation of
    ${_{\st{\psi_0}}}\Y_{\st{\A}} \in \Bimc(\B,\A)$.
    It is easy to check that the map $[\qX] \mapsto [\Y]$ of
    isomorphism classes is well defined.
\end{proof}

We remark that  $\cl_*$ is not surjective in general, see Section 
\ref{sec:imagecl}.

In the case where $\A = \B$ and $\star = \star'$, a simple computation
shows that if $\qX, \qX' \in \Picc(\qA)$ are $(\star,\star)$-bimodule
deformations of $\X,\X' \in \Picc(\A)$, then
$$
\cl_*([\qX \otimes_{\st{\qA}}\qX']) =  [\X\otimes_{\st{\A}}\X'].
$$ 
This observation implies the next result. 
\begin{lemma}
    \label{lem:clasmap}
    The classical-limit map (\ref{eq:clpic})
    \begin{equation}
        \label{eq:clasmap}
        \cl_*:\Pic(\qA) \longrightarrow \Pic(\A)
    \end{equation}
    is a group homomorphism. 
\end{lemma}

\begin{example}
    \label{example:trivdef}
    Consider the trivial deformation $\qA$ where the product is just
    the undeformed product of $\A$ extended $k[[\lambda]]$-bilinearily
    to $\qA$. Then 
    $$
    \Pic(\qA) = \PicA(\A) \ltimes \Aut(\qA)
    $$ 
    is larger than $\Pic(\A)$ due to the larger automorphism group, showing that
    (\ref{eq:clasmap}) is not injective in general.
\end{example}

Let us consider the group of self-equivalences of $\qA$,
\begin{equation}
   \label{eq:equivA}
    \Equiv(\qA) := \{T \in \Aut(\qA) \; | \; T = \id + O(\lambda)\}.
\end{equation} 
Note that, since $\A$ is commutative, $\Inaut(\qA) \subseteq
\Equiv(\qA)$ is a normal subgroup. We define the group of 
\emph{outer self-equivalences} of $\qA$ as
\begin{equation}
    \label{eq:outerA}
    \outequiv(\qA):= \frac{\Equiv(\qA)}{\Inaut(\qA)}.
\end{equation}
\begin{proposition}
    \label{prop:invimage}
    For any $[\X] \in \Bim(\B,\A)$ in the image of $\cl_*$,
    $\cl_*^{-1}([\X])$ is in one-to-one correspondence with the outer
    self-equivalences of $\qB$.
\end{proposition}
\begin{proof}
    Let $\qX$, $\qX' \in \Bimc(\qB,\qA)$, and let $[\X] =\cl_*([\qX])$
    and $[\X'] \in \cl_*([\qX'])$. If $\cl_*([\qX]) = \cl_*([\qX'])$,
    then, in particular, $\X_{\st{\A}} \cong \X'_{\st{\A}}$ as right
    $\A$-modules. So, by uniqueness of right-module deformations (see
    e.g. \cite{BuWa2000b}), $\qX_{\st{\qA}} \cong \qX'_{\st{\qA}}$ as
    right $\qA$-modules. Hence there exists a $\psi \in \Aut(\qB)$
    such that
    $$
    {_{\st{\psi}}}\qX_{\st{\qA}}
    \cong {_{\st{\qB}}}{\qX'}_{\st{\qA}}
    $$ 
    as $(\qB,\qA)$-bimodules.

    Note that 
    $\cl_*( [ {_{\st{\psi}}}\qX_{\st{\qA}}]) = [{_{\st{\psi_0}}}\X_{\st{\A}}]$, 
    where $\psi_0$ is such that $\psi = \psi_0 + O(\lambda)$.
    By Proposition~\ref{prop:transit}, 
    $_{\st{\psi_0}}\X_{\st{\A}} \cong { \AXA}$ if and only if $\psi_0 = \id$.
    So, again by Proposition~\ref{prop:transit}, 
    $$
    \cl_*^{-1}([\X]) 
    = \{ [_{\st{\psi}}\qX_{\st{\qA}}] \; |\; \psi \in \Equiv(\qB)\} 
    \cong \outequiv(\qB).
    $$
\end{proof}

For the group homomorphism $\cl_*:\Pic(\qA) \longrightarrow \Pic(\A)$
from (\ref{eq:clasmap}), we have $\ker(\cl_*) =
\cl_*^{-1}([{_{\st{\A}}}\A_{\st{\A}}])$.  As a consequence, we obtain
the following result.
\begin{corollary}
    \label{cor:ker}
    For the classical-limit map (\ref{eq:clasmap}),
    we have a group isomorphism $\ker(\cl_*) \cong \outequiv(\qA)$.
\end{corollary}

In fact, the group isomorphism is given by
$$
\ker(\cl_*)  \ni [{_{\st{\psi}}\qA}]
\mapsto [\psi^{-1}] \in \outequiv(\qA),
$$
see (\ref{eq:ldef}).

\subsection{Picard-group actions on deformations}
\label{sec:picaction}

In this section, we will present a description of the image of the map
$\cl_*$ in (\ref{eq:clasmap}).  Let $\A$ be a commutative unital
algebra.  We start with a slight extension of the discussion in
\cite{Bu2002a}.
\begin{lemma} 
    \label{lem:choices}
    Let $\star$ be a deformation of $\A$ , and let $\X, \Y$ be
    isomorphic objects in $\Picc(\A)$. Suppose
    $\qA'=(\A[[\lambda]],\star')$ and $\qA''=(\A[[\lambda]],\star'')$
    are deformations of $\A$ for which there exist a
    $(\star',\star)$-bimodule deformation $\qX \in \Bimc(\qA',\qA)$
    of $\X$, and a $(\star'',\star)$-bimodule deformation $\qY \in
    \Bimc(\qA'',\qA)$ of $\Y$. Then $[\star']=[\star'']$.
\end{lemma}
\begin{proof}
    Let $\psi_0:\X \longrightarrow \Y$ be a bimodule isomorphism.
    By uniqueness of right-module deformations up to equivalence
    \cite{BuWa2000b}, there exists a right $\qA$-module isomorphism
    $$
    \psi=\psi_0 + \sum_{r=1}^\infty \lambda^r \psi_r: \;
    \qX_{\st{\qA}} \longrightarrow \qY_{\st{\qA}}.
    $$
    Denote the left-module structure of $\qX$ over $\qA'$ (resp. $\qY$
    over $\qA'')$ by $\bullet'$ (resp. $\bullet'')$. For each 
    $a \in \qA''$, the map 
    $$
    \qX \ni x \;\mapsto\;  \psi^{-1}(a \bullet'' \psi(x)) \in \qX
    $$
    is right $\qA$-linear. Hence there is a unique element 
    $\phi(a) \in \qA'$ defined by
    \begin{equation}
        \label{eq:choice}
        \phi(a)\bullet' x 
        = \psi^{-1}(a \bullet'' \psi(x)), 
        \quad\textrm{for all}\; x \in \qX.
    \end{equation}
    It is easy to check that 
    $\phi=\phi_0+\sum_{r=1}^\infty\lambda^r \phi_r: 
    \qA'' \longrightarrow \qA'$ 
    is an algebra isomorphism. 
    Finally, notice that (\ref{eq:choice}) in zeroth order implies that
    $$
    \phi_0(a)x = \psi_0^{-1}(a\psi_0(x))
    = a x, 
    \quad\textrm{for all}\; a \in \A, \; x \in \X.
    $$
    Hence $\phi_0 = \id$, and $[\star'] = [\star'']$.
\end{proof}
\begin{theorem}
    \label{thm:action1}
    There exists a natural action 
    $\widetilde{\Phi}: \Pic(\A)\times \Def(\A) \longrightarrow \Def(\A)$
    so that two deformations $\star$ and $\star'$ of $\A$ are Morita
    equivalent if and only if $[\star]$ and  $[\star']$ lie in the
    same $\widetilde{\Phi}$-orbit.
\end{theorem}
\begin{proof}
    Let $\qA=(\A[[\lambda]],\star)$ be a deformation of $\A$, and 
    $\X \in \Picc(\A)$. We will first show that there is a deformation
    $\star'$ of $\A$ for which there exists a
    $(\star',\star)$-bimodule deformation of $\X$.

    Recall that, as a right $\A$-module, $\X$ can be identified with
    $P\A^n$ for some idempotent $P \in M_n(\A)$.  If $\qP=P +
    O(\lambda) \in M_n(\qA)$ is a (necessarily full) idempotent
    deforming $P$, then, by identifying $\qP\star \qA^n$ with
    $\X[[\lambda]]$, we obtain an induced module deformation of
    $\X_{\st{\A}}$ with respect to $\star$. We denote this deformation
    by $\qX_{\st{\qA}}$.

    On the other hand, we can identify $\A[[\lambda]]$ with
    $\End(\qX_{\st{\qA}})$ as $k[[\lambda]]$-modules, see
    e.g.~\cite{BuWa2000b}, to obtain an induced deformation
    $(\A[[\lambda]],\hat{\star})$ of $\A$, together with a left action
    $\hat{\bullet}$ of $(\A[[\lambda]],\hat{\star})$ on $\qX$. This
    turns $\qX$ into an $(\hat{\star},\star)$-bimodule. Note that, in
    general, the classical limit of $\qX$ is isomorphic to a bimodule
    of the form $_{\st{\psi}}\X$, for some $\psi \in \Aut(\A)$.
    Defining the deformation $\qA'=(\A[[\lambda]],\star')$, with
    $$
    \star' := (\psi^{-1})^*(\hat{\star}),
    $$
    we get a bimodule $_{\st{\qA'}}\qX_{\st{\qA}}$ with classical
    limit $\AXA$, thereby obtaining a $(\star',\star)$-bimodule
    deformation of $\AXA$. It is also clear by the construction that
    $_{\st{\qA'}}\qX_{\st{\qA}}$ is an equivalence bimodule.

    A simple application of Lemma~\ref{lem:choices} shows that
    $[\star']$ depends only upon the isomorphism class $[\X] \in
    \Pic(\A)$ and the equivalence class $[\star]\in \Def(\A)$.  So we
    have a well-defined map
    $$
    \widetilde{\Phi}: 
    \Pic(\A)\times \Def(\A) \longrightarrow \Def(\A), 
    \quad [\star] \mapsto 
    \widetilde{\Phi}_{\st{\X}}([\star]) = [\star'].
    $$

    To see that $\widetilde{\Phi}$ defines an action, let $\star$ be a
    deformation of $\A$. Note that for the identity $\A \in \Picc(\A)$,
    $\widetilde{\Phi}_{\st{\A}} = \id$, since $\qA$ is itself an
    $(\qA,\qA)$-equivalence bimodule deforming $\A$. Now let $\X, \Y
    \in \Picc(\A)$, and suppose $ \star' \in
    \widetilde{\Phi}_{\st{\X}}([\star]) $ and $\star'' \in
    \widetilde{\Phi}_{\st{\Y}}([\star'])$, in such a way that
    $$
    [\star''] = \widetilde{\Phi}_{\st{\Y}} 
    \widetilde{\Phi}_{\st{\X}}([\star]).
    $$
    Then there exists a $(\star',\star)$-bimodule deformation of $\X$
    and a $(\star'',\star')$-bimodule deformation of $\Y$, denoted by
    $\qX$ and $\qY$, respectively. But $\qY \otimes_{\st{\qA'}} \qX$
    is a $(\star'',\star)$-bimodule deformation of $\Y \tensorA \X$.
    So 
    $$
    [\star''] = \widetilde{\Phi}_{\st{\X\otimes_{\A}\Y}}([\star]).
    $$

    Finally, if $\star'$ and $\star$ are Morita equivalent, there
    exists a corresponding equivalence bimodule $\qX$. If 
    $[\X] = \cl_*([\qX])$, then $[\star'] 
    = \widetilde{\Phi}_{\st{\X}}([\star])$.
\end{proof}

We observe that, for a deformation $\star$ of $\A$, the image of the map
$\cl_*$ in (\ref{eq:clasmap}) is the isotropy group of $\widetilde{\Phi}$
at $[\star]$.

It turns out that the action $\widetilde{\Phi}$ can be written as a
combination of its restriction to the subgroup $\PicA(\A)$ and the
natural action of $\Aut(\A)$ on deformations of $\A$, see~\cite{Bu2002a}.
\begin{proposition}
    \label{prop:action2}
    There exists a natural action 
    $\Phi: \PicA(\A)\times \Def(\A) \longrightarrow 
    \Def(\A)$ such  that two deformations $\star, \star'$ of $\A$ are
    Morita equivalent if and only if there exists an automorphism
    $\psi$ of $\A$ so that $[\psi^*(\star')]$ and $[\star]$ lie in the
    same $\Phi$-orbit. 
\end{proposition}
\begin{proof}
    The action is just the one from Theorem~\ref{thm:action1},
    restricted to $\PicA(\A)$.
    So, if $[\psi^*(\star')]$ and $[\star]$ lie in the same
    $\Phi$-orbit, then $\psi^*(\star')$ and $\star$ are Morita
    equivalent. Since $\star'$ and $\psi^*(\star')$ are isomorphic
    deformations, it follows that $\star'$ is also Morita equivalent
    to $\star$.
    
    On the other hand, 
    if $\star'$ and $\star$ are Morita equivalent, then there exists a
    corresponding equivalence bimodule $\qX$,
    and $\cl_*([\qX]) = [{_{\st{\phi}}}\X] $  for some 
    $\phi \in \Aut(\A)$ and $\X \in \PiccA(\A)$. 
    Replacing $\star'$ by $\star''=(\phi^{-1})^*(\star')$,
    we obtain a $(\star'',\star)$-bimodule deformation of  $\X$.
    So $[\psi^*(\star')] = \Phi_{\st{\X}}([\star])$, where $\psi = \phi^{-1}$.
\end{proof}

As observed in \cite{Bu2002a}, the restricted action $\Phi$ has the nice feature 
that deformations in the same $\Phi$-orbit correspond to the same
Poisson bracket. So if we fix a Poisson bracket 
$\{\, , \, \}$ on $\A$, we have the following result.
\begin{corollary}
    \label{cor:action3}
    There is a natural action 
    $\Phi: \PicA(\A)\times \Def(\A,\{\,,\,\}) \longrightarrow 
    \Def(\A,\{\,,\,\})$ in such a way that two deformations 
    $\star, \star'$ of $(\A,\{\,,\,\})$
    are Morita equivalent if and only if there exists a Poisson
    automorphism $\psi$ of $\A$ so that $[\psi^*(\star')]$ and
    $[\star]$ lie in the same $\Phi$-orbit.
\end{corollary}

\begin{remark}
    \label{remark:thefuture}
    In fact, one can compute the semiclassical limit of the action
    $\Phi$, obtaining a first-order obstruction to Morita
    equivalence of formal deformations in terms of  
    algebraic Poisson cohomology \cite{Hueb90}; the arguments are similar
    to the case of star products
    \cite{Bu2002a}. A detailed analysis will be presented in a future work.
\end{remark}

\begin{remark}
    \label{rmk:staract}
    Consider again $\A = C^\infty(M)$, in which case $\PicA(\A) \cong
    \Pic(M) \cong H^2(M,\mathbb{Z})$, see Example \ref{ex:semidir}.
    For a star product $\star$ on $M$ and $[L] \in \Pic(M)$, one can
    show that there exists a $\star' \in \Phi_{\st{L}}([\star])$ which
    is again a star product on $M$, i.e.\ a deformation given by
    bidifferential cochains, see~\cite{BuWa2000b,Bu2002a}. So the
    action $\Phi$ in Proposition~\ref{prop:action2} restricts to an
    action
    \begin{equation}
        \label{eq:staract}
        \Phi: \Pic(M) \times \Def(M) \longrightarrow \Def(M),\;\; 
        [\star] \mapsto \Phi_{\st{L}}([\star]).
    \end{equation}
    For a fixed Poisson structure $\pi$ on $M$,
    two star products $\star$ and $\star'$ on $(M,\pi)$ are Morita
    equivalent if and only if there exists a Poisson diffeomorphism
    $\psi:M \longrightarrow M$ such that $[\psi^*(\star')]$ and
    $[\star]$ lie in the same $\Phi$-orbit, see \cite[Thm.~4.1]{Bu2002a}.
\end{remark}

Let us fix a deformation $\star$ of $\A$, corresponding to a Poisson
bracket $\{\, ,\, \}$. In order to describe the image of
$\cl_*:\Pic(\qA) \longrightarrow \Pic(\A)$, we may first consider the
map 
\begin{equation}
    \label{eq:resclass}
    \cl_*^r = \pr \circ \cl_*: \; \Pic(\qA) \longrightarrow \PicA(\A),
\end{equation}
where $\pr: \Pic(\A) \longrightarrow \PicA(\A)$ is the natural
projection, see Corollary~\ref{cor:semidir}.

By Proposition~\ref{prop:action2}, the image of $\cl_*^r$
is given by those $[\X] \in \PicA(\A)$ for which there exists a
$\star' \in \Phi_{\st{\X}}([\star])$ and $\psi \in \Aut(\A)$ such that
$[\star] = [\psi^*(\star')]$. Since $\star'$ and $\star$ correspond to
the same Poisson bracket $\{\, , \}$, it follows that $\psi$ must be a
Poisson map. As a result, we have the following corollary.
\begin{corollary}
    \label{cor:imagered}
    The image of the map $\cl_*^r$ in (\ref{eq:resclass}) is given by 
    $$
    \image(\cl_*^r)=\{ [\X] \in \PicA(\A) 
    \; | \; \exists \psi \in \Poiss(\A,\{\,,\,\}) \;\textrm{and}\;
    \star' \in \Phi_{\st{\X}}([\star])
    \mbox{ with } [\psi^*(\star')] = [\star]\},
    $$
    where $\Poiss(\A,\{\,,\,\})$ denotes the group of Poisson
    automorphisms of $(\A,\{\,,\,\})$.
\end{corollary}

In terms of the image of $\cl_*^r$, the image of $\cl_*$ can be
described as follows. For $[\X] \in \image(\cl_*^r)$, let
$$
P_{\st{\X}}:= \{ \psi \in \Poiss(\A,\{\,,\,\}) \;\textrm{for which} \;
\exists \, \star' \in \Phi_{\st{\X}}([\star]) \;\textrm{with}\;
[\psi^*(\star')] = [\star]\}.
$$
\begin{corollary}
    \label{cor:image}
    The image of $\cl_*$ is given by
    $$
    \image(\cl_*) = \{[{_{\st{\psi}}}\X]\; | \; [\X] \in
    \image(\cl_*^r) \;\textrm{and}\; 
    \psi \in P_{\st{\X}}\}.
    $$
\end{corollary}

Note that, since $[{_{\st{\psi}}}\X] \in \image(\cl_*)$ implies that
$\psi \in \Poiss(\A,\{\,,\,\})$ and $\Pic(\A) = \Aut(\A)\ltimes
\PicA(\A)$ (see  Corollary~\ref{cor:semidir}), it immediately follows
that $\cl_*$ is not onto in general.  We will come back to this
question in Section~\ref{sec:imagecl}.


\section{Bimodule deformations and their semiclassical limit}
\label{sec:bidef}

In this section, we identify the semiclassical limit of bimodule
deformations and discuss the corresponding deformation quantization
problem.

\subsection{Contravariant connections}
\label{sec:contra}

Let $\A$ be a commutative unital $k$-algebra, and let
$\{\,,\,\}$ be a Poisson bracket on $\A$. The next definition
follows the notion of contravariant derivatives in Poisson geometry
\cite{Rui2000,Hueb90,Vais94}.
\begin{definition}
    \label{def:contav}
    A \emph{contravariant connection} on $\X \in \PiccA(\A)$ is a
    $k$-bilinear map $D: \A \times \X \longrightarrow \X$ satisfying
    \begin{itemize}
    \item[i)] $D(ab,x)= D(a,x)b + D(b,x)a$,
    \item[ii)] $D(a,bx) = bD(a,x) + \{a,b\}x$.
    \end{itemize}
    We denote $D(a, \cdot) = D_a$. The \emph{curvature} of $D$ is the
    skew-symmetric bilinear map 
    $\curv(D): \A \times \A \longrightarrow \A$ defined by 
    $$
    \curv(D)(a,b) = D_aD_b - D_bD_a -D_{\{a,b\}}, \;\; a,b \in \A.
    $$
    Here, we are using the identification of $\End(\X)$ with
    $\A$. As usual, $D$ is called \emph{flat} if $\curv(D) = 0$.
\end{definition}
Note that if  $D$ and $D'$ are contravariant connections on 
$\X$, then their difference $D-D'$ is a derivation of $\A$, in the
sense that $D_ax - D'_a x = \alpha(a)x$ determines a derivation
$\alpha$ of $\A$. Conversely, if $\alpha : \A \longrightarrow \A$
is a derivation and $D$ is a contravariant connection, 
so is $D + \alpha$.

\begin{example}[\cite{Vais94,Rui2000}]
    Let $(M,\pi)$ be a Poisson manifold, $L$ a line bundle over $M$
    and $\nabla$ any (ordinary) connection on $L$. For $f \in C^\infty(M)$, let
    $X_f$ denote its Hamiltonian vector field.
    Then 
    $$
    D: C^\infty(M) \times \Gamma^\infty(L) \longrightarrow \Gamma^\infty(L)
    $$
    defined by $D(f,s) = \nabla_{X_f} s$ is a contravariant
    connection. When $\pi$ is symplectic, any contravariant connection
    arises in this way.
\end{example}

\begin{example}\label{ex:d}
For the $(\A,\A)$-bimodule $\A$, with respect to left and right multiplication,
we can define the canonical contravariant connection $d$ by
$$
d_a(b)=\{a,b\}, \;\; a,b \in \A.
$$
This connection is flat as a consequence of the Jacobi identity of $\{\, , \,\}$.
\end{example}

Let $\X, \X' \in \PiccA(\A)$ be equipped with contravariant
connections $D$ and $D'$, respectively.
\begin{definition}
    \label{def:isoconn}
    An \emph{isomorphism of bimodules preserving contravariant
      connections} is a bimodule isomorphism $T: \X \longrightarrow
    \X'$ satisfying
    $$
    T(D_a x) = D_a' Tx,
    $$ 
    for all $a \in \A$ and $x\in \X$.
\end{definition}
Note that isomorphic contravariant connections have the same
curvature. We denote the set of contravariant connections on $\X$
with curvature $\tau$ by $\Connc(\X,\tau)$; the set 
of isomorphism classes of elements in
$\Connc(\X,\tau)$ is denoted by $\Conn(\X,\tau)$.

If $\alpha$ is a derivation of $\A$, we call it a 
\emph{Poisson derivation} if 
$$
\alpha(\{a,b\}) = \{ \alpha(a), b \} + \{ a, \alpha(b) \}.
$$
We denote the set of derivations and Poisson derivations of $\A$ by 
$\Der(\A)$ and $\PDer(\A)$, respectively.
\begin{lemma}
    \label{lem:samecurv}
    Let $D$ be a contravariant connection on $\X \in \PiccA(\A)$, and
    let $\alpha$ be a derivation of $\A$. Then $\curv(D+\alpha) =
    \curv(D)$ if and only if $\alpha \in \PDer(\A)$. In particular,
    $\Connc(\X, \tau)$ is an affine space over $k$ modeled on
    $\PDer(\A)$.
\end{lemma}
\begin{proof}
    If we  expand the equation
    $$
    (D+\alpha)_a(D+\alpha)_b(x)-(D+\alpha)_b(D+ \alpha)_a(x)
    -(D+\alpha)_{\{a,b\}}(x) =\tau,
    $$
    it is simple to check that $\curv(D) = \tau$ if and only if  $\alpha$
    is a Poisson derivation of $\A$.
\end{proof}
\begin{proposition}
    \label{prop:equivconn}
    Let $\X \in \PiccA(\A)$, $D \in \Connc(\X,\tau)$, 
    and $\alpha \in \Der(\A)$. Then $(\X,D + \alpha)$
    is isomorphic to $(\X,D)$ if and only if $(\A,d + \alpha)$ is
    isomorphic to $(\A,d)$.
\end{proposition}
\begin{proof}
    Suppose there exists a bimodule automorphism 
    $\phi: \A \longrightarrow \A$ satisfying
    $$
    \phi \, d_a = (d_a + \alpha(a)) \, \phi,
    $$
    for all $a \in \A$.
    For the bimodule with contravariant connection
    $(\X,D)$, consider the tensor products $(\X\otimes \A, D\otimes d)$
    and $(\X \otimes \A, D\otimes (d+ \alpha))$, where the tensor
    product of two connections is defined in the natural way:
    $$
    D_1\otimes D_2 (x_1 \otimes x_2) = D_1(x_1)\otimes x_2 + x_1 \otimes D_2(x_2).
    $$
    One can check that the natural identification $\X\otimes \A
    \longrightarrow \X$, $x \otimes a \mapsto xa$ induces isomorphisms
    $$
    (\X\tensorA \A, D\otimes d) \stackrel{\sim}{\longrightarrow} (\X,D),
    $$
    and
    $$
    (\X \tensorA \A, D \otimes (d + \alpha))
    \stackrel{\sim}{\longrightarrow} (\X, D + \alpha).
    $$
    Finally, note that $\phi$ gives rise to a map $\X\tensorA \A
    \longrightarrow \X \tensorA \A$ through $x \otimes a \mapsto x
    \otimes \phi(a)$, which establishes an isomorphism between
    $D\otimes d$ and $D \otimes (d + \alpha)$. Hence there is an
    induced isomorphism
    $$
    (\X, D) \stackrel{\sim}{\longrightarrow} (\X, D + \alpha).
    $$
    
    For the converse, suppose that there is an isomorphism
    $$
    \psi : (\X, D) \longrightarrow (\X, D+\alpha).
    $$
    Define the dual module with contravariant connection $(\X^*,D^*)$, where
    $D^*$ and $D$ are related in the usual way:
    $$
    d_a\SP{x^*,x} = \SP{D^*_a x^*, x} + \SP{x^*,D_a x},
    $$
    where $a \in \A$, $x \in \X$, $x^* \in \X^* =\Hom(\X,\A) $ and
    $\SP{\; , \;}$ is the pairing between $\X$ and $\X^*$.  Then
    $\psi$ induces an isomorphism
    $$
    (\X \tensorA \X^*, D \otimes D^*)
    \stackrel{\sim}{\longrightarrow} (\X\tensorA \X^*,(D+\alpha)\otimes D^*).
    $$
    Finally, notice that the isomorphism $\X \tensorA \X^*
    \longrightarrow \A$ given by $x \otimes x^* \mapsto \SP{x^*,x}$
    induces isomorphisms $(\X \tensorA \X^*, D \otimes D^*)
    \stackrel{\sim}{\longrightarrow} (\A,d)$ and $(\X\tensorA
    \X^*,(D+\alpha)\otimes D^*)\stackrel{\sim}{\longrightarrow} (\A,
    d+\alpha)$. Therefore $(\A,d)$ and $(\A,d+\alpha)$ are isomorphic.
\end{proof}
\begin{corollary}
    \label{cor:classflat}
    Isomorphism classes of contravariant connections on $\X$ with
    fixed curvature $\tau$ are classified by isomorphism classes of
    flat contravariant connections on $\A$. 
\end{corollary}

We have the following characterization of Poisson derivations yielding
isomorphic connections.
\begin{lemma}
    \label{lem:equivconn}
    A contravariant connection $D + \alpha$ on $\X$ is isomorphic to $D$
    if and only if there exists an invertible element $u \in \A$ such that
    $\alpha = u^{-1}\{\cdot, u\}$.
\end{lemma}
\begin{proof}
    By Corollary~\ref{cor:classflat}, it suffices to consider the case
    $\X = \A$, $D=d$. Any  bimodule automorphism of $\A$ is given by
    multiplication by an invertible element $u \in \A$. Note that
    $ u (d + \alpha)_a b = d_a b $ for all $a,b \in \A$, if and only
    if 
    $$
    u(\{a,b\} + \alpha(a)b) = \{a, ub\} = \{a,b\}u + \{a, u\}b,
    $$ 
    which is equivalent to  $\alpha(a) = \{a,u\}/u $.
\end{proof}
\begin{definition}
    \label{def:intder}
    Let $u \in \A$ be invertible. We call  the derivations of the form
    $\alpha = u^{-1}\{u, \cdot\}$ \emph{integral derivations}.
    We denote the set of integral derivations of $\A$ by
    $\IDer(\A)$.
\end{definition}
Note that if $uv=1$, then $u^{-1}\{u, \cdot\} = v^{-1}\{\cdot,v\}$, so
we have two alternative definitions of integral derivations.  Note
also that integral derivations form an abelian group: the addition
$\alpha+\alpha'$ in $\IDer(\A)$ corresponds to the product $uu'$.
Heuristically speaking, $\alpha$ corresponds to the inner Poisson
derivation with the (not necessarily existing) logarithm of $u$.

Since isomorphic connections have the same curvature, it follows from
Lemma~\ref{lem:samecurv} that $\IDer(\A) \subseteq \PDer(\A)$. Of
course, this can also be checked directly.
\begin{corollary}\label{cor:classconn}
    The set $\Conn(\A,0)$ is an affine space over $\mathbb{Z}$ modeled on
    $\PDer(\A)\big/\IDer(\A)$ with canonical origin $[d]$. Hence
    the map 
    $$
    \Conn(\A, 0) \longrightarrow
    \PDer(\A)/\IDer(\A),\;\;\;
    [D] \mapsto [D - d]
    $$
    is a bijection. By Corollary~\ref{cor:classflat}, $\Conn(\X,\tau)$ is
    also an affine space over $\mathbb{Z}$ modeled on
    $\PDer(\A)\big/\IDer(\A)$, but with no canonical origin in
    general. 
\end{corollary}
\begin{example}
    \label{ex:sym}
    Let $(M,\omega)$ be a symplectic manifold and
    $\A = C^\infty(M)$, equipped with the induced Poisson bracket.
    We have
    $$
    \IDer(\A) \cong \{ f^{-1} df  \;|\; f:M \to \mathbb{C}^* \}
    \subseteq Z^1(M,\mathbb{C)},
    $$
    where $Z^1(M,\mathbb{C})$ denotes the space of closed
    complex-valued $1$-forms on $M$.

    Recall that on a contractible open set $\mathcal{O}_\alpha
    \subseteq M$, any smooth $f:M \longrightarrow \mathbb{C}^*$ can be
    written as $f = \eu^{2\pi\im c_\alpha}$, for $c_{\alpha} \in
    C^\infty(\mathcal{O}_\alpha)$.  Hence $f^{-1}df = 2\pi\im
    dc_\alpha$ and $c_\alpha - c_\beta \in \mathbb{Z}$ on overlaps
    $\mathcal{O}_\alpha \cap \mathcal{O}_\beta$.  It follows that
    $$
    \pr(\IDer(\A)) = 2\pi\im \HdR^1(M,\mathbb{Z}),
    $$
    where $\pr: Z^1(M,\mathbb{C)} \longrightarrow
    \HdR^1(M,\mathbb{C})$ denotes the usual projection and
    $\HdR^1(M,\mathbb{Z})$ is the image of $H^1(M,\mathbb{Z})$ under
    the natural map
    $$
     H^1(M,\mathbb{Z}) \longrightarrow \HdR^1(M,\mathbb{C}).
    $$
    This motivates our terminology of integral derivations.  Note that
    $$
    \PDer(\A)/{\IDer(\A)} \cong 
    Z^1(M,\mathbb{C})/\{ f^{-1} df \;|\; f:M \to \mathbb{C}^* \},
    $$
    and since $\{ f^{-1} df \;|\; f:M \to \mathbb{C}^* \}$ contains all
    exact $1$-forms, there is a well-defined surjective map
    $$
    \HdR^1(M,\mathbb{C}) \to 
    \frac{Z^1(M,\mathbb{C)}}{\{ f^{-1} df \;|\; f:M \to \mathbb{C}^* \}}.
    $$
    The kernel of this map is just 
    $\pr(\IDer(\A)) = 2 \pi \im \HdR^1(M,\mathbb{Z})$.
    So, as abelian groups,
    \begin{equation}
        \label{quots}
        \frac{\PDer(\A)}{{\IDer(\A)}} \cong 
        \frac{\HdR^1(M,\mathbb{C})}{2\pi\im \HdR^1(M,\mathbb{Z})}.
    \end{equation}
    Corollary~\ref{cor:classconn} implies the well-known fact that
    isomorphism classes of lines bundles with connections of fixed
    curvature are classified by the quotient (\ref{quots}).
\end{example}
\begin{example}
    \label{ex:poiss}
    Let $(M,\pi)$ be a general Poisson manifold.  In this case,
    integral Poisson derivations are Poisson vector fields of the form
    $$
    f^{-1}\{\cdot , f\} = f^{-1}d_{\pi}f,
    $$
    for $f \in C^\infty(M)$, where $d_{\pi} = [\pi,\,\cdot\,]$ and
    $[\;,\;]$ is the Schouten bracket \cite{Vais94,SilWein99}.
    Consider the map
    $$
    \tilde{\pi}:\Omega^1(M) \to \chi^1(M), \quad
    \tilde{\pi}(\beta)=\pi(\cdot,\beta).
    $$
    Since $d_{\pi}f = \tilde{\pi}(df)$, there is an induced map in
    cohomology $\pi^*:\HdR^1(M,\mathbb{C}) \longrightarrow
    H^1_{{\pi}}(M,\mathbb{C})$, where
    $H^1_{{\pi}}(M,\mathbb{C})$ denotes the first Poisson
    cohomology of $M$, see e.g.~\cite{SilWein99}.  For $f \in
    C^\infty(M)$, we have
    $$
    f^{-1}d_{\pi}f = \tilde{\pi}(f^{-1}df),
    $$ 
    and hence
    $$
    \{[f^{-1}d_{\pi}f]_{\pi},\;\; f:M\to \mathbb{C}^*\} =
    \{ \pi^*[f^{-1}df],\;\; 
    f:M\to \mathbb{C}^*\}  
    = 2\pi\im H^1_{\pi}(M,\mathbb{Z}),
    $$
    where $H^1_\pi(M, \mathbb{Z}) := \pi^*\HdR^1(M,\mathbb{Z})$.
    Analogously to the previous example, we conclude that
    \begin{equation}
        \label{eq:quotp}
        \frac{\PDer(\A)}{\IDer(\A)} \cong 
        \frac{H^1_{\pi}(M,\mathbb{C})}{2\pi\im H^1_{\pi}(M,\mathbb{Z})}
    \end{equation}
    as abelian groups. Corollary~\ref{cor:classconn} asserts the fact
    that (\ref{eq:quotp}) classifies isomorphism classes of line
    bundles over $M$ with contravariant connections of fixed
    curvature.
\end{example}

\subsection{The semiclassical limit of bimodule deformations}

Let $\X \in \PiccA(\A)$. Let $\star$ be a deformation of
$(\A,\{\,,\,\})$, and let $\star' \in \Phi_{\st{\X}}([\star])$.
We write these deformations as
\begin{equation}
    \label{eq:star}
    \star = \sum_{r=0}^\infty \lambda^r C_r 
    \quad\textrm{and}\quad
    \star' = \sum_{r=0}^\infty \lambda^r C_r'.
\end{equation}
\begin{definition}
    \label{def:tau}
    We define the skew-symmetric bilinear map $\tau: \A \times \A \longrightarrow \A$,
    corresponding to $\star', \star$, by
    $$
    \tau(a_1,a_2) = (C_2-C_2')(a_1,a_2) - (C_2-C_2')(a_2,a_1),\;\;\; a_1,a_2 \in \A.
    $$
\end{definition}
Since $\A$ is commutative and $\star$ and $\star'$ deform the same
Poisson bracket, $\tau$ can also be written as
\begin{equation}
    \label{eq:tau2}
    \tau(a_1,a_2) 
    = \left.\left(
            \frac{1}{\lambda^2}
            \left([a_1,a_2]_{\star} - [a_1,a_2]_{\star'}\right)
        \right)
    \right|_{\lambda = 0}
\end{equation}
By Proposition~\ref{prop:action2} the algebras
$\qA=(\A[[\lambda]],\star)$ and $\qA'=(\A[[\lambda]],\star')$ are
Morita equivalent. Let
$$
\qX=(\X[[\lambda]],\bullet',\bullet) \in \Bimc(\qA',\qA)
$$ 
be a $(\star',\star)$-bimodule deformation of $\X$, and let us write
$$
\bullet' = \sum_{r=0}^\infty \lambda^r R_r' 
\quad\textrm{and}\quad
\bullet = \sum_{r=0}^\infty \lambda^r R_r.
$$ 
\begin{proposition}
    \label{prop:semiclas}
    Suppose  $C_1=C_1'$. Then the  map
    $$
    D = R_1 - R_1': \A \times \X \longrightarrow \X
    $$ 
    is a contravariant connection, and $\curv(D) = \tau$.
\end{proposition}
\begin{remark}\label{rem:conv}
    In the case of star products, in order to be consistent with the
    convention (\ref{eq:convention}), we define $D = (R_1 -
    R_1')/\im$; in this case $\curv(D) = -\tau$.
\end{remark}
\begin{proof}
    The proof is analogous to the proofs of \cite[Prop.~4.3 and
    Thm~5.3]{Bu2002a} in the case of star products. We consider the
    equations relating $\star'$, $\star$, $\bullet'$ and $\bullet$:
    \begin{eqnarray}
        (a_1\star' a_2) \bullet' x 
        & = & a_1 \bullet'(a_2 \bullet' x),  \label{eq:rel1}\\
        x \bullet(a_1 \star a_2) 
        & = & (x \bullet a_1) \bullet a_2, \label{eq:rel2}\\
        (a_1\bullet' x)\bullet a_2 
        & = & a_1\bullet'(x\bullet a_2), \label{eq:rel3}
    \end{eqnarray}
    for $a_1, a_2 \in \A$, $x \in \X$. A suitable combination of
    equations (\ref{eq:rel1}), (\ref{eq:rel2}) and (\ref{eq:rel3}) in
    order $\lambda$ shows that $D=R_1-R_1'$ is a contravariant
    connection. A combination of these equations in order $\lambda^2$
    implies that $\curv(D)=\tau$, see \cite{Bu2002a} for details.
\end{proof}

In particular, $\curv(D)$ depends  only upon $\star'$ and $\star$, and
not on any particular $(\star',\star)$-bimodule deformation of $\X$.

Let $\Defc(\X,\star',\star)$ be the set of
$(\star',\star)$-bimodule deformations of $\X$; the set of isomorphism classes
of elements in $\Defc(\X,\star',\star)$ is denoted by $\Def(\X,\star',\star)$.

Using
Proposition~\ref{prop:semiclas}, we define the semiclassical limit
map
\begin{equation}
    \label{eq:S}
    S: \Defc(\X,\star',\star) \longrightarrow \Connc(\X,\tau), \quad
    (\X[[\lambda]],\bullet',\bullet) \mapsto D=R_1- R_1'.
\end{equation}

\begin{example}\label{ex:Sd}
Consider the $(\A,\A)$-bimodule $\A$ with respect to left and
right multiplication. In this case, the contravariant connection
$S(\qA)$ associated to the $(\star,\star)$-bimodule deformation of
$\A$ given by $\qA$ is just $d$ of Example \ref{ex:d},
\begin{equation}
    \label{eq:d}
    d_a(b) = C_1(a,b) - C_1(b,a) = \{a,b\}. 
\end{equation}
\end{example}

\begin{proposition}
    \label{prop:equivsemi}
    Let $\qX_1, \qX_2$ be $(\star',\star)$-bimodule deformations of
    $\X_1, \X_2 \in \PiccA(\A)$, respectively. Let $D_1=S(\qX_1)$ and
    $D_2=S(\qX_2)$ be the induced contravariant connections on $\X_1$
    and $\X_2$. If
    $$
    T=\sum_{r=0}^\infty \lambda^r T_r : \qX_1 \longrightarrow \qX_2
    $$
    is a bimodule isomorphism, then $T_0:\X_1 \longrightarrow \X_2$
    is an isomorphism preserving contravariant connections. In
    particular, if $T$ is an equivalence, $D_1 = D_2$.
\end{proposition}
\begin{proof}
    Let us denote the left and right actions on $\qX_i$ by 
    $$
    \bullet_i' = 
    \sum_{r=0}^\infty \lambda^r R_r^{i'} 
    \quad\textrm{and}\quad
    \bullet_i =\sum_{r=0}^\infty \lambda^r R_r^{i},
    \;\; i=1,2.
    $$
    From $T(a \bullet_1' x)= a \bullet_2' T(x)$ and
    $T(x \bullet_1 a) = T(x) \bullet_2 a$, we get, in order $\lambda$,
    \begin{eqnarray*}
        T_1(a x) + T_0(R^{1'}_1(x,a))
        & = & R_1^{2'}(a,T_0(x)) + a T_1(x), \\
        T_1(x a) +T_0(R_1^1(x,a))
        & = & R_1^2(T_0(x),a) + T_1(x)a.
    \end{eqnarray*}
    Subtracting the two equations and recalling that we are assuming
    that the classical left and right actions are the same, we get
    $$
    T_0({D_1}(a,x)) = T_0(R_1^1(x,a)- R^{1'}_1(x,a))=
    R_1^2(T_0(x),a)- R_1^{2'}(a,T_0(x))= {D_2}(a,T_0(x)).
    $$
\end{proof}

As a result of Proposition~\ref{prop:equivsemi}, the semiclassical
limit map $S$ (\ref{eq:S}) induces a map of isomorphism classes
\begin{equation}
    \label{eq:Stilde}
    \widetilde{S}: \Def(\X,\star',\star) \longrightarrow \Conn(\X,\tau).
\end{equation}
\begin{remark}
    Let $\star$ be a deformation of $\A$ and $\X \in \PiccA(\A)$.  Let
    $D$ be a contravariant connection on $\X$.  A natural question is
    whether one can find $\star' \in \Phi_{\st{\X}}([\star])$, with
    $C_1'=C_1$, and a $(\star',\star)$-bimodule deformation of $\X$
    with $D$ as the semiclassical limit. For classes of deformations
    for which Hochschild $2$-cochains are cohomologous to their
    skew-symmetric parts (e.g, for star products), one can always find
    $\star' \in \Phi_{\st{\X}}([\star])$ with $C_1'=C_1$ since $\Phi$
    preserves Poisson brackets. In this case, a simple computation
    shows that one can choose $\star'$ so that there exists a
    $(\star',\star)$-bimodule deformation $\qX$ of $\X$ with $S(\qX) =
    D$, see \cite{Bu2002b}.
\end{remark}

\subsection{Deformation quantization of bimodules}
\label{sec:defquantbim}

Let $\X \in \PiccA(\A)$, $\star$ be a deformation of $(\A,\{\,,\,\})$
and $\star' \in \Phi_{\st{\X}}([\star])$.  Following the notation in
(\ref{eq:star}), we assume that $C_1=C_1'$, and consider the
semi-classical limit map  
$$
S:\Defc(\X,\star',\star) \longrightarrow \Connc(\X,\tau)
$$
as in (\ref{eq:S}).
\begin{definition}
    \label{def:quantconn}
    A contravariant connection $D$ on $\X$ is called \emph{quantizable}
    if $D$ is in the image of $S$.
\end{definition}
Analogously to deformation of algebras, we call an element 
$\qX \in S^{-1}(D)$ a $(\star',\star)$-bimodule deformation of $\X$
\emph{in the direction of} $D$.

In this section, we will address the following two deformation
quantization problems: 
\begin{enumerate}
\item First, whether any contravariant connection on $\X$ with
    curvature $\tau$ (determined by $\star'$ and $\star$) is
    quantizable (i.e, whether $S$ is onto);
\item Second, how to classify equivalence classes of
    $(\star',\star)$-bimodule deformations of $\X$ in the direction of
    a fixed $D \in \image(S) $, i.e., equivalence classes of elements
    in $S^{-1}(D)$.
\end{enumerate}

Let $\star$ be a deformation of $\A$, and consider the map
\begin{equation}
    \label{eq:s}
    s: \Equiv(\qA) \longrightarrow \PDer(\A), \;\; \;
    \id + \lambda T_1 + O(\lambda^2) \mapsto T_1,
\end{equation}
which clearly is a group homomorphism. Let $u = u_0 + O(\lambda) \in
\qA$ be invertible.  A simple computation shows that $s(\Ad(u)) =
\{u_0, \cdot\}/{u_0}$.  Since any invertible $u_0 \in \A$ is also
invertible in $\qA$, it follows that
\begin{equation}
    \label{eq:incl}
    s(\Inaut(\qA)) = \IDer(\A).
\end{equation}
\begin{definition}
    \label{def:qder}
    A Poisson derivation $\alpha \in \PDer(\A)$ is called
    \emph{quantizable} if there exists a $\star$-derivation $\deriv$
    with $\deriv = \alpha + O(\lambda)$.
\end{definition}
\begin{remark}
    \label{rem:exp}
    Recall that any $\star$-self-equivalence $T$ is of the form
    $\eu^{\lambda \deriv}$ for a unique $\star$-derivation $\deriv$
    (see e.g.~\cite[Lem.~5]{BuWa2002}).  So $\alpha \in \PDer(\A)$ is
    quantizable if and only if it lies in the image of $s$.  By
    (\ref{eq:incl}), any integral derivation of $\A$ is quantizable.
\end{remark}
\begin{theorem}\label{thm:defcond}
    Let $\X \in \PiccA(\A)$, $\star$ be a deformation of
    $(\A,\{\,,\,\})$ and $\star' \in \Phi_{\st{\X}}([\star])$.  With the
    notation of (\ref{eq:star}), let us assume that $C_1=C_1'$. The
    semiclassical limit map 
    $S:\Defc(\X,\star',\star) \longrightarrow \Connc(\X,\tau)$ 
    is onto if and only if 
    $s: \Equiv(\qA') \longrightarrow \PDer(\A)$ is onto.
\end{theorem}
\begin{proof}
    For $\qX \in \Defc(\X,\star',\star)$, let $D = S(\qX)$.  For each
    $T= \id + \lambda T_1 + O(\lambda) \in \Equiv(\qA')$, a simple
    computation shows that $S({_T}\qX) = D - T_1$. This shows that if
    any Poisson derivation is quantizable, we can choose $T$
    conveniently and, by Lemma \ref{lem:samecurv}, 
    any contravariant connection on $\X$ with
    curvature $\tau$ is quantizable.

    Now suppose $S$ is onto.  Let $\alpha$ be a Poisson derivation.
    Then there exists a deformation $\qY$ of $\X$ so that $S(\qY) = D
    + \alpha$. On the other hand, $\qY \cong {_T}\qX$ for some $T \in
    \Equiv(\qA')$.  Since $S({_T}\qX) = D - T_1$, it follows from
    Proposition~\ref{prop:equivsemi} that $D + \alpha$ and $D - T_1$
    are isomorphic connections on $\X$. By Lemma~\ref{lem:equivconn},
    this is the case if and only if
    $$
    \alpha = - T_1 + u^{-1}\{\cdot, u \},
    $$
    for some invertible $u \in \A$.  But since both $T_1$ and
    $u^{-1}\{\cdot, u \}$ are quantizable, so is $\alpha$.
\end{proof}

Keeping the notation of Theorem~\ref{thm:defcond}, we have the
following result concerning classification of bimodule deformations.
\begin{theorem}
    \label{thm:classif}
    Let $D \in S(\Defc(\X,\star',\star))$. The set of equivalence
    classes of bimodule deformations in $S^{-1}(\X,D)$ is in
    one-to-one correspondence with the set
    \begin{equation}
        \label{eq:classif}
        \frac{\{ T \in \Equiv(\qA') \;|\; T_1=0\}}
        {\{ \Ad(u) \in \Inequiv(\qA') \;|\; u = 1 + O(\lambda) \in \qA'\}}.
    \end{equation}
\end{theorem}
\begin{proof}
    If $\qX$ is a deformation of $(\X,D)$, then, for $T=\id + \lambda
    T_1 + O(\lambda^2) \in \Equiv(\qA')$, ${_T}\qX$ deforms $(\X, D -
    T_1)$. So ${_T}\qX$ deforms $(X,D)$ if and only if $T_1=0$.
    
    The bimodules $\qX$ and ${_{\st{T}}\qX}$ are isomorphic if and
    only if $T = \Ad(u)$, for some invertible $u= u_0 + O(\lambda) \in
    \qA'$. The induced bimodule isomorphism $\qX \longrightarrow
    {_{\st{T}}\qX}$ is given by
    $$
    x \mapsto u \bullet' x, \;\;\; x \in \qX.
    $$
    So this isomorphism is an equivalence if and only if $u_0 = 1$.    
\end{proof}


\section{Inner automorphisms of star product algebras}
\label{sec:auto}

In this section we briefly discuss the structure of the automorphism group
of a star product $\star$ on a
Poisson manifold $(M, \pi)$. As usual, we let $\A=C^\infty(M)$ and
$\qA = (C^\infty(M)[[\lambda]], \star)$.

\subsection{The star exponential}
\label{sec:starexp}

We shall need the star exponential \cite{BFFLS78} as a technical tool
and will recall its main properties below.

For a given $H \in \qA$ the \emph{star exponential} $\Exp(H)$ is
defined as the $t=1$ value of the unique solution to the differential
equation
\begin{equation}
    \label{eq:starexpdiff}
    \frac{d}{dt} \Exp(tH) = H \star \Exp(tH)
\end{equation}
with initial condition $\Exp(0) = 1$. It is well-known that
(\ref{eq:starexpdiff}) has indeed a solution with the following
properties, see e.g.~\cite[App.~A]{BuWa2002}.
\begin{lemma}
    \label{lemma:Exp}
    Let $H = \sum_{r=0}^\infty \lambda^r H_r \in
    C^\infty(M)[[\lambda]]$. Then (\ref{eq:starexpdiff}) has a unique
    solution $\Exp(tH)$ of the form
    \begin{equation}
        \label{eq:ExpExpl}
        \Exp(H) = \sum_{r=0}^\infty \lambda^r \Exp(H)_r
        \quad\textrm{with}\quad
        \Exp(H)_r \in C^\infty(M),
    \end{equation}
    where $\Exp(H)_0 = \eu^{H_0}$ and 
    $\Exp(H)_{r+1} = \eu^{H_0} H_{r+1} + E_r$ with $E_r$ depending on
    $H_0, \ldots, H_r$ only. The following additional properties hold.
    \begin{enumerate}
    \item $H \star \Exp(H) = \Exp(H) \star H$ and 
        $\Exp(tH) \star \Exp(sH) = \Exp((t+s)H)$.
    \item Let $M$ be connected. Then $\Exp(H) = 1$ if and only if $H$
        is constant and equal to $2\pi\im m$, for some
        $m \in \mathbb{Z}$.
    \item For all $f \in C^\infty(M)[[\lambda]]$ one has 
        $\Exp(H) \star f \star \Exp(-H) = \eu^{\ad(H)} (f)$, where $\ad(H)f
        = [H,f] = H\star f - f \star H$ is the $\star$-commutator.
    \item $[f,g] = 0$ if and only if $[\Exp(f), \Exp(g)] = 0$ if and
        only if $[f, \Exp(g)] = 0$. In this case,
        $\Exp(f) \star \Exp(g) = \Exp(f+g)$.
    \end{enumerate}
\end{lemma}
\begin{lemma}
    \label{lemma:LocalLog}
    Let $\mathcal{O} \subseteq M$ be an open contractible subset and
    let $U \in C^\infty(\mathcal{O})[[\lambda]]$ be invertible. Then
    there exists a function $H \in C^\infty(\mathcal{O})[[\lambda]]$
    such that $\Exp(H) = U$ and any two such functions differ by a
    constant $2\pi\im m$, with $m \in \mathbb{Z}$.
\end{lemma}

\subsection{Locally inner derivations}
\label{sec:locinnder}

Let $\deriv: \qA \to \qA$ be a $\mathbb{C}[[\lambda]]$-linear derivation of
$\star$. We call $\deriv$ \emph{locally inner} if, on each contractible
open subset $\mathcal{O} \subseteq M$, there exists a function $H \in
C^\infty(\mathcal{O})[[\lambda]]$ such that
\begin{equation}
    \label{eq:locinnder}
    \deriv \big|_{C^\infty(\mathcal{O})[[\lambda]]} = \ad(H).
\end{equation}
The $\mathbb{C}[[\lambda]]$-module of locally inner derivations is
denoted by $\LocInnDer(\qA)$.
\begin{definition}
    \label{definition:defdeltaA}
    Let $A \in Z^1(M,\mathbb{C})[[\lambda]]$ be a closed $1$-form.  We
    define a global $\star$-derivation $\delta_A$ by
    \begin{equation}
        \label{eq:deltaADef}
        \delta_A\big|_{C^\infty(\mathcal{O})[[\lambda]]} = \ad(H),
    \end{equation}
    where $\mathcal{O}$ is a contractible open subset of $M$ and 
    $H$ is a local function with $dH = A$ on $\mathcal{O}$.
\end{definition}
It is clear that $\delta_A$ is well-defined. It is a locally inner
derivation by construction and gives rise to a
$\mathbb{C}[[\lambda]]$-linear map
\begin{equation}
    \label{eq:delta}
    \delta: Z^1(M,\mathbb{C})[[\lambda]] \longrightarrow 
    \LocInnDer(\qA),
    \quad A \mapsto \delta_A.
\end{equation}

In the symplectic case, the map $\delta$ is known to be a bijection
between closed $1$-forms and the star product derivations, see
e.g.~\cite{GR99}. In the Poisson case, $\delta$ is generally neither
injective nor surjective. Moreover, the (non-)bijectivity of $\delta$
depends on $\star$ itself and not only on the Poisson bracket as the
following example illustrates.
\begin{example}
    \label{example:SymplBadGuy}
    Consider a star product $\star'$ on a symplectic manifold, and let
    $\star$ be defined by replacing $\lambda$ by $\lambda^2$ in
    $\star'$.  Then $\star$ is a deformation of the zero Poisson
    bracket and $\delta$ is still bijective. Note that, in this case,
    the Poisson center and the center of $\star$ do not coincide.  On
    the other hand, for the trivial star product corresponding to the
    zero Poisson bracket, any vector field is a derivation but
    $\delta_A$ is identically zero.
\end{example}

The situation becomes  easier when $\star$ has trivial center.
\begin{lemma}
    \label{lemma:TrivCenter}
    Suppose the center of $(C^\infty(M)[[\lambda]],\star)$ is trivial.
    Then $\delta$ is injective and surjective onto the set of locally
    inner derivations.
\end{lemma}
\begin{proof}
    Let $\{\mathcal{O}_{\alpha}\}$ be a cover of $M$ by contractible
    open sets.  Given a locally inner derivation $\deriv$ with $\deriv
    = \ad(H_\alpha)$ on $\mathcal{O}_\alpha$ and $\deriv =
    \ad(H_\beta)$ on $\mathcal{O}_\beta$, we see that $H_\alpha -
    H_\beta$ is a central function on $\mathcal{O}_\alpha \cap
    \mathcal{O}_\beta$ hence locally constant. Thus $A = dH_\alpha$
    defines a \emph{global} $1$-form and $\deriv = \delta_A$.  On the
    other hand, if $\delta_A = 0$ then $\ad(H) = 0$. So $H$ is
    constant and $0 = dH = A$.
\end{proof}

We finally observe the following simple result.
\begin{lemma}
    \label{lemma:locinnerideal}
    Let $A \in Z^1(M,\mathbb{C})[[\lambda]]$, and let $\deriv$
    be an arbitrary $\star$-derivation. Then $[\deriv, \delta_A]$ is inner.
\end{lemma}
\begin{proof}
    Obvious, as $\deriv H_\alpha = \deriv H_\beta$ for $A = dH_\alpha = dH_\beta$
    on $\mathcal{O}_\alpha \cap \mathcal{O}_\beta$.
\end{proof}

\subsection{Inner automorphisms of star products}
\label{sec:innerauto}

We will now characterize the inner automorphism of the star-product
algebra $\qA = (C^\infty(M)[[\lambda]],\star)$.
\begin{theorem}
    \label{prop:edeltaA}
    \begin{enumerate}
    \item Let $u=\sum_{r=0}^\infty\lambda^r u_r \in \qA$ be invertible
        and let $T=\Ad(u)$. Then there exists a closed $1$-form $A \in
        Z^1(M, \mathbb{C})[[\lambda]]$ with $T=\eu^{\delta_A}$.
        Moreover, $A_0$ is $2 \pi \im$-integral and the $A_r$ are
        exact for $r \geq 1$.
    \item The derivations $\LogInnAut(\qA)$ of $\qA$ which
        exponentiate to inner automorphisms form an abelian group
        under addition and we have the following inclusions of groups
        \begin{equation}
            \label{eq:groupsinothergroups}
            \InnDer(\qA)
            \subseteq \LogInnAut(\qA)
            \subseteq \delta(Z^1(M, \mathbb{C})[[\lambda]])
            \subseteq \LocInnDer(\qA).
        \end{equation}
    \item The map
        \begin{equation}
            \label{eq:IntegralHoneToLogINnAut}
            2 \pi \im \HdR^1 (M, \mathbb{Z}) 
            \ni [A] \mapsto [\delta_{A}] \in
            \frac{\LogInnAut(\qA)}{\InnDer(\qA)}        
        \end{equation}
        is a surjective group homomorphism. If the center of $\star$ is
        trivial, it is a bijection.
    \end{enumerate}
\end{theorem}
\begin{proof}
    For the first part, recall that $T = \eu^{\lambda \deriv}$ with a
    unique $\star$-derivation $\deriv$, see Remark~\ref{rem:exp}.  Let
    $\{\mathcal{O}_\alpha\}$ be a cover of $M$ by contractible open
    sets.  By Lemma~\ref{lemma:LocalLog}, there exist local functions
    $H_\alpha \in C^\infty(\mathcal{O}_\alpha)[[\lambda]]$ with $u =
    \Exp(H_\alpha)$ on $\mathcal{O}_\alpha$. Hence $T =
    \eu^{\ad(H_\alpha)}$ on $\mathcal{O}_\alpha$.  Since the functions
    $H_\alpha$ are unique up to constants, the $1$-form $A =
    dH_\alpha$ is well-defined globally. Thus $T = \eu^{\delta_A}$.
    Finally, the $H_\alpha$ differ by constants in $2\pi\im\mathbb{Z}$
    whence $A_0$ is $2\pi\im$-integral, and the $A_r$ are exact for $r
    \ge 1$. The second part follows from the first part. The third
    part is a consequence of (\ref{eq:groupsinothergroups}) and
    Lemma~\ref{lemma:TrivCenter}.
\end{proof}

Let $\cente_{\pi}(\A)$ denote the Poisson center of $(\A,\{\;,\;\})$, and let
$\cente(\qA)$ denote the $\star$-center of $\qA$. Then we define the
following property of a star product $\star$:
\begin{enumerate}
\item[C0)] There exists a linear map
    $\rho_0:\cente_\pi(\A) \longrightarrow \cente(\qA)$ with
    $\rho_0(f) = f + O(\lambda)$.
\end{enumerate}
Note that, in this case, $\rho_0$ actually extends to a
$\mathbb{C}[[\lambda]]$-linear bijection between
$\cente_\pi(\A)[[\lambda]]$ and the center $\cente(\qA)$. Hence the
center $\cente(\qA)$ is a deformation of the sub-space $\cente_\pi(\A)$
in the sense of \cite[Def.~30]{bordemann.herbig.waldmann:2000}.  In
Section~\ref{sec:linedef} we will discuss examples of star products
for which such a $\rho_0$ exists.
\begin{theorem}
    \label{theorem:InnerAutoStar}
    Let $T=\Ad(u)$ be an inner automorphism of $\qA$ of the form $T =
    \id + O(\lambda^2)$ and suppose that $\star$ satisfies C0).  Then
    there exists $u'=1 + O(\lambda)$ such that $T = \Ad(u') =
    \eu^{\ad(H)}$ for some global function $H \in \lambda
    C^\infty(M)[[\lambda]]$.
\end{theorem}
\begin{proof}
    Let us write $T(f) = u \star f \star v$. Evaluating $u \star f
    \star v$ in first order, which has to vanish according to our
    assumptions, we get the condition
    $$
    u_0fv_1 + u_1 f v_0 + C_1(u_0,f)v_0 + u_0C_1(f,v_0) = 0.
    $$
    Since $T^{-1}$ is also of the form $\id + O(\lambda^2)$, we
    obtain the same equation with the roles of $u$ and $v$ exchanged.
    Taking their difference, we get
    $$
    v_0\{u_0,f\} = -u_0\{f,v_0\}.
    $$
    Using $u_0v_0=1$, we conclude that
    $$
    \{u_0,f\}=0  \; \mbox{ for all } \; f \in C^\infty(M),
    $$
    and therefore $u_0 \in \cente_{\pi}(\A)$. So $\rho_0(u_0) = u_0
    + O(\lambda) \in \cente(\qA)$. Since $u_0$ is invertible, so is
    $\rho_0(u_0)$, and we can define
    $$
    u'= \rho_0(u_0)^{-1}\star u.
    $$
    Note that $v' = \rho_0(u_0)\star v$ is the $\star$-inverse of
    $u'$, and clearly $T = \Ad(u)=\Ad(u')$ as $\rho_0(u_0)$ is
    central.  Since in this case the formal series of the star
    logarithm trivially converges in the $\lambda$-adic topology, we
    can find $H\in \lambda C^\infty(M)[[\lambda]]$ with $u' =
    \Exp(H)$.  Hence $T(f) = \Exp(H) \star f \star \Exp(-H) =
    \eu^{\ad(H)}(f)$.
\end{proof}
\begin{corollary}\label{cor:outerder}
    If $\star$ satisfies C0) then there exists a bijection 
    \begin{equation}
        \label{eq:OutDerBijection}
        \frac{\{ T \in \Equiv(\qA)\;|\; T_1=0\}}
        {\{\Ad(u) \in \Inequiv(\qA) \;|\; u =1 + O(\lambda)\}}
        \longrightarrow \mathrm{qOutDer}(\qA),
    \end{equation}
 where $\mathrm{qOutDer}(\qA)$ is the quotient $(\lambda\Der(\qA))/\InnDer(\qA)$.
\end{corollary}
\begin{proof}
There is a one-to-one correspondence between 
self-equivalences of $\qA$ and $\star$-derivations given by
$$
T= \eu^{\lambda \deriv} \mapsto \deriv.
$$
Under this bijection, the set $\{ T \in \Equiv(\qA),\; T_1=0\}$ is
mapped bijectively onto derivations of the form $\lambda \deriv$, and,
by Theorem~\ref{theorem:InnerAutoStar}, the set $\{\Ad(u) \in
\Inequiv(\qA), \; u =1 + O(\lambda)\}$ corresponds to inner
derivations. So the result follows.
\end{proof}

As the lefthand side of (\ref{eq:OutDerBijection}) is a group, one may
ask which group structure on $\mathrm{qOutDer}(\qA)$ is obtained  from
this bijection. This is illustrated in the following example.
\begin{example}
    \label{example:trivdefout}
    Let $\star$ be the trivial deformation, i.e. the undeformed
    product. Then there are no non-trivial inner automorphisms and the
    equivalence transformations with $T_1 = 0$ are of the form
    $\eu^{\lambda\mathcal{L}_X}$ with some vector field $X \in
    \lambda\Gamma^\infty(TM)[[\lambda]]$. Thus 
    $\eu^{\lambda\mathcal{L}_X} \eu^{\lambda\mathcal{L}_Y}
    = \eu^{\lambda\mathcal{L}_{H(X,Y)}}$, where 
    $$
    H(X,Y) = X + Y + \frac{1}{2}[X,Y] + \cdots
    $$
    is the Baker-Campbell-Hausdorff series, which is well-defined
    as $X, Y$ start in order $\lambda$. Here $\mathcal{L}_{X}$ is the
    Lie derivative of $X$. Thus $\mathrm{qOutDer}(\qA)$ inherits the
    non-abelian group structure
    $$
    X \cdot Y = \frac{1}{\lambda}H(\lambda X, \lambda Y).
    $$
\end{example}

\section{Bimodule deformations of line bundles over Poisson manifolds}
\label{sec:linedef}

We will now specialize the discussion of deformation quantization of
invertible bimodules to the context of line bundles over a Poisson
manifold $(M, \pi)$. We will consider star products $\star$ on $M$
satisfying property C0) and the following analogous property for
derivations:
\begin{enumerate}
\item[C1)] There exists a linear map $\rho_1: \PDer(\A)
    \longrightarrow \Der(\qA)$, with $\rho_1(X) = \mathcal{L}_X +
    O(\lambda)$, such that, for Hamiltonian vector fields, 
    $\rho_1(X_H) = \frac{\im}{\lambda}\ad(H)$. 
\end{enumerate}
The map $\rho_1$ extends to a $\mathbb{C}[[\lambda]]$-linear bijection
between $\PDer(\A)[[\lambda]]$ and $\Der(\qA)$ yielding a deformation
of the subspace $\PDer(\A)$.

We call $\star$-derivations of the form $\frac{\im}{\lambda}\ad(H)$
\emph{quasi-inner}; note that the quotient of $\Der(\qA)$ by
quasi-inner derivations can be naturally identified with $\mathrm{qOutDer}(\qA)$. 
Since we assume that $\rho_1$ maps Hamiltonian
vector fields into quasi-inner $\star$-derivations, a simple induction
shows the following lemma.
\begin{lemma}\label{lem:outerpoisson}
    Let $\star$ be a star product on $(M,\pi)$ satisfying C1). Then
    $\mathrm{qOutDer}(\qA)$ is in one-to-one correspondence with
    $H_{{\pi}}^1(M,\mathbb{C})[[\lambda]]$.
\end{lemma}

Before we move on, let us remark that one can always find star
products on $(M,\pi)$ satisfying conditions C0) and C1).
\begin{example}
    Let $(M,\omega)$ be a symplectic manifold, and let $\star$ be any
    star product on $M$.  In this case, the Poisson center of $\A$ is
    trivial and one can take $\rho_0$ to be the canonical embedding.  As
    mentioned in Section~\ref{sec:locinnder}, the map $\delta$ in
    (\ref{eq:delta}) is a bijection in the symplectic case.  The
    correspondence between symplectic vector fields and closed one
    forms,
    $$
    X \mapsto \contract_X\omega=\omega(X, \,\cdot\, ),
    $$
    induces a one-to-one correspondence 
    $$
    \hat{\delta}:\chi_\omega(M)[[\lambda]] \longrightarrow \Der(\qA), \;\; 
    X \mapsto \frac{\im}{\lambda}\delta_A, 
    \quad\textrm{with}\quad
    A = \contract_X\omega,
    $$
    where $\chi_\omega(M)$ denotes the space of symplectic vector
    fields.  Using the convention in (\ref{eq:convention}), a simple
    computation shows that
    $$
    \hat{\delta}_X = \mathcal{L}_{X} + O(\lambda).
    $$
    It is easy to check that $\hat{\delta}$ maps Hamiltonian vector
    fields on $M$ onto quasi-inner $\star$-derivations.  So one can
    take $\hat{\delta}$ as $\rho_1$.  As a consequence, we have a
    group isomorphism between $\mathrm{qOutDer}(\qA)$ and elements
    in $\HdR^1(M,\mathbb{C})[[\lambda]]$, see e.g.~\cite{GR99}.

    In order to describe the group structure on
    $\HdR^1(M,\mathbb{C})[[\lambda]]$ induced by
    (\ref{eq:OutDerBijection}), let $A,A'$ be closed $1$-forms.
    From Lemma~\ref{lemma:locinnerideal} and the
    Baker-Campbell-Hausdorff theorem, it follows that
    $\eu^{\delta_A} \eu^{\delta_{A'}} = \eu^{\delta_{A+A'}}$ modulo
    some inner automorphism. Thus the induced group structure on
    $\mathrm{qOutDer}(\qA)$ is just the abelian one, quite in contrast
    with Example~\ref{example:trivdefout}.
\end{example}
\begin{example}
    Let $(M,\pi)$ be an arbitrary Poisson manifold, and let
    $\mathcal{U}$ be a formality on $M$ \cite{Kon97}, i.e., an
    $L_\infty$-quasi-isomorphism from the differential graded Lie
    algebra of multivector fields on $M$ (with zero differential and
    Schouten bracket) into the differential graded Lie algebra of
    multidifferential operators on $M$ (with Hochschild differential
    and Gerstenhaber bracket).  It is well-known that $\mathcal{U}$ is
    completely determined by its Taylor coefficients $\mathcal{U}_p$,
    $p\geq 1$. Kontsevich showed in \cite{Kon97} that formalities
    exist on arbitrary Poisson manifolds and that one can take
    $\mathcal{U}_1$ to be the natural embedding of vector fields into
    differential operators.  Given $\mathcal{U}$, one can define a
    star product on $(M,\pi)$ by
    \begin{equation}\label{eq:kontsp}
        f \star g 
        := fg + \sum_{r=1}^\infty \frac{\lambda^r}{r!}
        \mathcal{U}_r(\pi, \ldots, \pi)(f \otimes g),
        \quad f,g \in C^\infty(M). 
    \end{equation}
    We call a star product associated to a formality as in
    (\ref{eq:kontsp}) a \emph{Kontsevich star product}.  More
    generally, Kontsevich has shown that any star product on $(M,\pi)$
    is equivalent to one obtained as in (\ref{eq:kontsp}) with $\pi$
    replaced by a \emph{formal Poisson structure} $\pi_{\lambda} = \pi +
    O(\lambda)$.
    
    As a consequence of the formality equations, one can show that
    Kontsevich star products satisfy C0) and C1) (see e.g.
    \cite{CFT01,CFT02,JSW2000,Gam2001}) for
    $$
    \rho_0(f) = f + \sum_{r=1}^\infty
    \frac{\lambda^r}{r!}\mathcal{U}_{r+1}(f, \pi,  \ldots, \pi),
    $$
    and
    $$
    \rho_1(X) = \mathcal{L}_X  + \sum_{r=1}^\infty
    \frac{\lambda^r}{r!}\mathcal{U}_{r+1}(X, \pi, \ldots, \pi).
    $$
\end{example}    
\begin{example}
    A more geometric way of obtaining global star products on
    $(M,\pi)$ is described in \cite{CFT01,CFT02}, where the authors
    use the formality on $\mathbb{R}^d$ combined with a Fedosov-type
    procedure \cite{Fed94}. It is shown in \cite{CFT01} that these star products
    can also be constructed satisfying C0) and C1). We will refer to
    these star products as \emph{CFT star products}.
\end{example}

\begin{theorem}\label{thm:injsympl}
    Let $\star$ be a star product on $(M,\pi)$ satisfying C0) and C1)
    (in particular, symplectic, CFT and Kontsevich star products).
    Let $L \to M$ be a complex line bundle, and let $\star' \in
    \Phi_{\st{L}}([\star])$.  Let $D$ be a contravariant connection on
    $L$ with curvature $-\tau$ determined by $\star$ and $\star'$.  
    Then there exists a
    $(\star',\star)$-bimodule deformation $\qX$ of
    $\X=\Gamma^\infty(L)$, with $S(\qX) =D$. Moreover, equivalence
    classes of bimodule deformations in $S^{-1}(D)$ are in one-to-one
    correspondence with $H_{{\pi}}^1(M,\mathbb{C})[[\lambda]]$.
\end{theorem}
\begin{proof}
    For the existence part, note that property C0) guarantees that the
    map $s$ (\ref{eq:s}) is onto. So the result follows from
    Theorem~\ref{thm:defcond}.  As for the second assertion, note
    that, since $\star$ satisfies C0), Corollary~\ref{cor:outerder}
    holds. Since property C1) is satisfied, we can apply
    Lemma~\ref{lem:outerpoisson}, and the result follows.    
\end{proof}

\section{The map $\cl_*$ for star-product algebras}
\label{sec:kerimagecl}

\subsection{The kernel of $\cl_*$}
\label{sec:kercl}

Let $(M,\pi)$ be a Poisson manifold and $\star$ be a star product on
$M$ satisfying C0) and C1).  As usual, we let $\A = C^\infty(M)$ and
$\qA=(C^\infty(M)[[\lambda]],\star)$.  Using $\rho_1$, we consider the
map
\begin{equation}
    \label{eq:equivp}
    \Equiv(\qA) \longrightarrow \PDer(\A)\times {\{T\in \Equiv(\qA),\; T_1=0\}},
    \;\;\;
    T \mapsto (T_1, \eu^{-\lambda \rho_1(T_1)}\circ T),
\end{equation}
which is easily seen to be a bijection.  
Let us  assume that the map (\ref{eq:equivp}) restricts to a
bijection
\begin{equation}
    \label{eq:inequivp}
    \Inequiv(\qA) \longrightarrow  \IDer(A) \times {\{T\in \Inequiv(\qA),\; T_1=0\}},
\end{equation}
which is always the case for symplectic star products if
one takes $\rho_1=\hat{\delta}$.

By Theorem \ref{theorem:InnerAutoStar}, (\ref{eq:inequivp}) actually
establishes a bijection
\begin{equation}\label{eq:inequivp3}
    \Inequiv(\qA) \longrightarrow  
    \IDer(\A) \times {\{T=\Ad(u) \in \Inequiv(\qA),\; u=1 +O(\lambda)\}}.
\end{equation}
Note that (\ref{eq:equivp}) and (\ref{eq:inequivp}) imply the existence of a a bijection
$$
\outequiv(\qA) \longrightarrow 
\frac{\PDer(\A)}{\IDer(\A)}\times
\frac{\{T\in \Equiv(\qA),\; T_1=0\}}{{\{T=\Ad(u),\; u_0 =1\}}}
\cong \frac{H^1_{{\pi}}(M,\mathbb{C})}{2\pi \im H^1_{{\pi}}(M,\mathbb{Z})}+
\lambda H^1_{{\pi}}(M,\mathbb{C})[[\lambda]],
$$
where the last identification follows from (\ref{eq:quotp}),
Corollary \ref{cor:outerder} and Lemma \ref{lem:outerpoisson}.  As a
consequence of Corollary \ref{cor:ker}, we can state the following
theorem.
\begin{theorem}
    For such a star product $\star$ 
    (e.g., any symplectic star product), there
    is a one-to-one correspondence
    \begin{equation}\label{eq:thmker}
    \ker(\cl_*) \stackrel{\sim}{\longrightarrow}
    \frac{H^1_{{\pi}}(M,\mathbb{C})}
    {2 \pi \im H^1_{{\pi}}(M,\mathbb{Z})}+ 
    \lambda H^1_{{\pi}}(M,\mathbb{C})[[\lambda]].
    \end{equation}
    In particular, $\cl_*$ is injective if and only if
    $H^1_{{\pi}}(M,\mathbb{C})=0$.
\end{theorem}
We remark that the identification (\ref{eq:thmker}) is a canonical group
isomorphism for symplectic star products. In this case, $\cl_*$ is injective if
and only if $\HdR^1(M,\mathbb{C})=0$.

However, in general (\ref{eq:thmker}) is not a group morphism for the
canonical abelian group structure on the right-hand side: If we consider
Example~\ref{example:trivdefout}, we see that the group structure on the left-hand 
side of (\ref{eq:thmker}) is  non-abelian.

\subsection{The image of $\cl_*$ for symplectic star products}
\label{sec:imagecl}

Let $(M,\omega)$ be a symplectic manifold. As shown in
\cite{BCG97,Deligne95,NT95b,WX98}, there exists a bijection
\begin{equation}
    \label{eq:starcharacclass}
   c : \Def(M,\omega) \longrightarrow \frac{1}{i \lambda}[\omega] +
   \HdR^2(M)[[\lambda]],
\end{equation}
characterizing the moduli space $\Def(M,\omega)$ in terms of the de Rham
cohomology of $M$. For a star product $\star$ on $M$, $c(\star)$
is called its \emph{characteristic class}.

Let $\star$ be a star product on $M$ with characteristic class
$$
    c(\star) 
    = \frac{1}{\im \lambda}[\omega] 
    + \sum_{r=0}^\infty \lambda^r [\omega_r].
$$
It is well known that the action of the group of
symplectomorphisms of $M$ on star products, in terms of
characteristic classes, is given by \cite{GR99}
\begin{equation}
    \label{eq:natural}
    c(\psi^*(\star))= \psi^*c(\star) = 
    \frac{1}{\im \lambda}[\omega] + \sum_{r=0}^\infty \lambda^r [\psi^* \omega_r].
\end{equation}
The simple description of this action is the main reason why we
restrict ourselves to symplectic star products in this section.

The main result in \cite{BuWa2002} asserts that, if $L \to M$ is a
complex line bundle and $\star' \in \Phi_{\st{L}}([\star])$ (see
Remark~\ref{rmk:staract}), then
\begin{equation}
    \label{eq:actioncharac}
    c(\star') = c(\star) + 2 \pi i c_1(L),
\end{equation}
where $c_1(L) \in H^2_{\st{dR}}(M,\mathbb{C})$ is the image under the canonical map
\begin{equation}\label{eq:e}
e: H^2(M,\mathbb{Z}) \longrightarrow \HdR^2(M,\mathbb{C})
\end{equation}
of the Chern class of $L$ in $H^2(M,\mathbb{Z})$.
Combining these results with Corollary \ref{cor:imagered}, we obtain
the following:
\begin{corollary}
    \label{cor:imclr}
    Let $\star$ be a star product on $(M,\omega)$, with $c(\star) =
    (1/{i \lambda})[\omega] + \sum_{r=0}^\infty \lambda^r [\omega_r]$.
    Then
    $$
    \image(\cl_*^r) = \{ l \in H^2(M,\mathbb{Z}) \; | \; \exists
    \psi \in \Sympl(M) \;\textrm{so that}\; e(l) =
    \psi^*[\omega_0]-[\omega_0] \mbox{ and } \psi^*[\omega_r] =
    [\omega_r],\, r \geq 0 \}.
    $$
\end{corollary}

For $l \in \image(\cl_*^r)$, let
$$
P_l=\{\psi \in \Sympl(M) \;|\; e(l) = \psi^*[\omega_0]-[\omega_0]
\;\textrm{and}\; \psi^*[\omega_r] = [\omega_r],\, r \geq 0 \}.
$$
Recall from Example~\ref{ex:semidir} that $\Pic(C^\infty(M)) =
\Diff(M) \ltimes \Pic(M)$, with action given by pull-back.
\begin{corollary}
   \label{cor:imclstar}
   The image of $\cl_*$ in $\Pic(C^\infty(M)) =
   \Diff(M)\ltimes\Pic(M)$ is given by  
   $$
   \image(\cl_*) = \{(\psi,l) \in \Diff(M) \ltimes \Pic(M) \; | \;
    \psi \in P_l \quad\textrm{and}\quad l \in \image(\cl_*^r) \}.
   $$
\end{corollary}

Let $\Tor(H^2(M,\mathbb{Z}))$ denote the subgroup of torsion elements
of $H^2(M,\mathbb{Z})$.  We recall that $\Tor(H^2(M,\mathbb{Z})) =
\ker(e)$. This subgroup can also be described as isomorphism classes
of line bundles over $M$ admitting a flat connection.  The following
result is a consequence of the Corollaries~\ref{cor:imclr} and
\ref{cor:imclstar}.
\begin{corollary}
   \label{cor:tor}
   Let $\star$ be a star product on $M$ such that $[\psi^*(\star)]
   =[\star]$ for all $\psi \in \Sympl(M)$. Then
   $$
   \image(\cl_*^r) = \Tor(H^2(M,\mathbb{Z}))\; \mbox{ and }\;
   \image(\cl_*) = \Sympl(M)\ltimes \Tor(H^2(M,\mathbb{Z})),
   $$
   with action given by pull-back.  Moreover, by
   Theorem~\ref{thm:injsympl}, if $\HdR^1(M,\mathbb{C})=0$, then
   $\Pic(\qA) \cong \Sympl(M)\ltimes \Tor(H^2(M,\mathbb{Z}))$.
\end{corollary}

Let us discuss a few examples.
\begin{example}
    Let us consider $M=\mathbb{R}^{2n}$, equipped with its canonical
    symplectic form.  In this case, $\PicA(\A) \cong
    H^2(\mathbb{R}^{2n},\mathbb{Z})=0$ and
    $$
    \Pic(\A) \cong \Diff(\mathbb{R}^{2n}).
    $$
    By Corollary~\ref{cor:tor}, for any symplectic star product
    $\star$ on $\mathbb{R}^{2n}$,
    $$
    \Pic(\qA) \cong \Sympl(\mathbb{R}^{2n}).
    $$
\end{example}    
\begin{example}
    More generally, for any symplectic manifold $(M,\omega)$, we can
    always choose a star product $\star$ satisfying $[\psi^*(\star)] =
    [\star]$ for all $\psi \in \Sympl(M)$: for example, one can take
    $\star$ with the trivial class
    $$
    c(\star) = \frac{1}{\im\lambda}[\omega].
    $$
    In this case, $\image(\cl_*^r)= \Tor(H^2(M,\mathbb{Z}))$ and
    $\image(\cl_*) = \Sympl(M)\ltimes \Tor(H^2(M,\mathbb{Z}))$.
\end{example}
\begin{example}
    Suppose $(M,\omega)$ is a symplectic manifold such that
    $\HdR^2(M,\mathbb{C})$ is one dimensional and $[\omega] \neq 0$.
    Then, for any $\psi \in \Sympl(M)$, the pull-back map
    $$
    \psi^*:\HdR^2(M,\mathbb{C}) \longrightarrow
    \HdR^2(M,\mathbb{C})
    $$
    must be the identity. As a result of (\ref{eq:natural}), for
    any star product $\star$,
    $$
    \psi^*([\star]) = [\star]\; \mbox{ for all }\, \psi \in
    \Sympl(M).
    $$
    Hence, by Corollary~\ref{cor:imclr}, $\image(\cl_*^r) =
    \Tor(H^2(M,\mathbb{Z}))$ and $\image(\cl_*) = \Sympl(M)\ltimes
    \Tor(H^2(M,\mathbb{Z}))$.

    In particular, if $M$ is a Riemann surface (equipped its area
    form) or $\mathbb{C}P^n$ (with the Fubini-Study symplectic form),
    then for any star product $\star$ on $M$ we have
    $$
    \image(\cl_*^r)=0 \quad\textrm{and}\quad \image(\cl_*) \cong
    \Sympl(M).
    $$
    In the case $M = \mathbb{C}P^n$, since
    $\HdR^1(\mathbb{C}P^n,\mathbb{C})=0$, we have
    $$
    \Pic(\qA) \cong \Sympl(\mathbb{C}P^n)
    $$
    for any star product.
\end{example}
\begin{example}
    In general, however, $\image(\cl_*^r)$ may contain non-torsion
    elements:  Let $\Sigma$ be a Riemann surface, and let $M =
    T^*\Sigma$, equipped with its canonical symplectic form. Recall
    that the natural projection $\pr:M\longrightarrow \Sigma$ induces
    an isomorphism
    $$
    \pr^*:\HdR^2(\Sigma,\mathbb{C}) \longrightarrow
    \HdR^2(M,\mathbb{C})
    $$
    For $f \in \Diff(\Sigma)$, let $f^\sharp \in \Sympl(M)$ be its
    cotangent lift.  We have the following commutative diagram:
    \[
    \begin{CD}
        \HdR^2(\Sigma,\mathbb{C}) @> {f^*} >>
        \HdR^2(\Sigma,\mathbb{C})  \\
        @V{\pr^*}VV  @VV{\pr^*}V\\
        \HdR^2(M,\mathbb{C}) @>{(f^\sharp)^*}>> \HdR^2(M,\mathbb{C}).
    \end{CD}
    \]
    Let
    $$
    \Diff(\Sigma)^*:=\{f^*:\HdR^2(\Sigma,\mathbb{C})\longrightarrow
    \HdR^2(\Sigma,\mathbb{C})\;,\; f \in \Diff(\Sigma)\}.
    $$
    Fixing an orientation in $\Sigma$, we can identify
    $\Diff(\Sigma)^*$ with $\mathbb{Z}_2 = \{\id,-\id\}$. Therefore
    $$
    \Sympl(M)^*:= \{\psi^*:\HdR^2(M,\mathbb{C})\longrightarrow
    \HdR^2(M,\mathbb{C})\;,\; \psi \in \Sympl(M)\} = \{\id,-\id\}.
    $$
    Let us pick $[\omega_0] \in \HdR^2(M,\mathbb{Z})$, and let $\star$
    be a star product on $M$ with
    $$
    c(\star)= \frac{1}{\im\lambda}[\omega] + [\omega_0] = [\omega_0],
    $$
    since $[\omega]=0$.  Then, by Corollary~\ref{cor:imclr},
    $$
    \image(\cl^r_*)= \{0, -2[\omega_0]\}.
    $$
    From Corollary~\ref{cor:imclstar}, it follows that
    $$
    \image(\cl_*)= A \cup B,
    $$
    where
    \begin{eqnarray*}
        A&=& \{(\psi,0) \in \Sympl(M)\ltimes H^2(M,\mathbb{Z})\,|\, \psi^*=\id: 
        \HdR^2(M,\mathbb{C}) \longrightarrow \HdR^2(M,\mathbb{C}) \}, \\
        B& =& \{(\phi,-2[\omega_0]) \in \Sympl(M)\ltimes H^2(M,\mathbb{Z})\,|\, \phi^*=-\id:
        \HdR^2(M,\mathbb{C}) \longrightarrow \HdR^2(M,\mathbb{C}) \}.
    \end{eqnarray*}
    It is simple to check that $A\cdot A \subseteq A$, $B \cdot A \subseteq B$,
    $A \cdot B \subseteq B$, and $B\cdot B \subseteq A$, and this shows
    directly that $\image(\cl_*)$ is a subgroup of $\Pic(\A)=\Sympl(M)\ltimes H^2(M,\mathbb{Z})$.
\end{example}

More complicated examples illustrating the nontriviality of
$\image(\cl^r_*)$ are obtained, e.g., by considering $M=T^*({T}^3)$ or
$M={T}^4$.  It seems hard to describe exactly how big the image of
$\cl_*^r$ can be in general.  We have, however, the following result.
\begin{theorem}
    Let $\star$ be a star product on $(M,\omega)$, and suppose
    $H^2_{\st{dR}}(M,\mathbb{C})$ is finite dimensional. Then the map
    $\cl_*^r$ is onto if and only if $\Pic(M)=H^2(M,\mathbb{Z})$ only
    contains torsion elements.
\end{theorem}
\begin{proof}
    It follows from Corollary~\ref{cor:imclr} that
    $\Tor(H^2(M,\mathbb{Z})) \subseteq \image(\cl_*^r)$, and this
    shows the ``if'' direction.
    
    For the converse, it suffices to show that
    $$
    H^2(M,\mathbb{Z}) \nsubseteq \{ l \in H^2(M,\mathbb{Z}) \; | \;
    \exists \psi \in \Sympl(M) \mbox{ so that } e(l) =
    \psi^*[\omega_0]-[\omega_0] \}.
    $$
    Fixing a basis of $\HdR^2(M,\mathbb{Z})$, we can identify
    $\HdR^2(M,\mathbb{C})$ with $\mathbb{C}^m$, and
    $\HdR^2(M,\mathbb{Z})$ with the lattice $\mathbb{Z}^m$.  Let us
    assume $H^2(M,\mathbb{Z})$ has non-torsion elements, so that
    $m\geq 1$.
    
    For $\psi\in \Sympl(M)$, $\psi^*$ must preserve the integer
    lattice, and therefore $\Sympl(M)^*\subseteq \GL(m,\mathbb{Z})$.
    So it suffices to show that, for any fixed $v \in \mathbb{C}^m$,
    $$
    \mathbb{Z}^m \nsubseteq \{ Av - v \;, \mbox{ for } A \in
    \GL(m,\mathbb{Z}) \}.
    $$
    
    Suppose $v=(v_1, \ldots, v_m) \in \mathbb{C}^m$ is such that, for
    any $l \in \mathbb{Z}^m$, there exists $A \in \GL(m,\mathbb{Z})$
    with $l = Av - v$. We will first show that, in this case, there
    exists a vector $r_0 \in \mathbb{Q}^m$ with the same property as
    $v$.
    
    Let $E = \mathrm{span}_{\mathbb{Q}}\{v_1, \ldots, v_m \}$. It
    follows from our assumptions that $1 \in E$. So we can find a
    basis for the $\mathbb{Q}$-vector space $E$ of the form $\{1, x_1,
    \ldots, x_s\}$. Using this basis, we write
    $$
    v = r_0 + r_1 x_1 + \cdots + r_s x_s,
    $$
    with $r_i \in \mathbb{Q}^m$.  Now suppose $(A-\id)v =l$, for $A
    \in \GL(m,\mathbb{Z})$ and $l \in \mathbb{Z}^m$.  Then
    $$
    (A-\id)r_0 + x_1(A - \id)r_1 + \cdots + x_s(A - \id)r_s = l,
    $$
    which implies that $x_1(A - \id)r_1 + \cdots + x_s(A - \id)r_s
    \in \mathbb{Q}^m$. But since $\{1, x_1, \ldots, x_s\}$ is a basis
    of $E$ over $\mathbb{Q}$, it must follow that
    $$
    x_1(A - \id)r_1 + \cdots + x_s(A - \id)r_s = 0.
    $$
    So if $l \in \mathbb{Z}^m$ and $(A-\id)v = l$, then $(A-\id)r_0
    =l$.
    
    We will now check that we cannot have $\mathbb{Z}^m$ contained in
    the set ${Av-v, \; A \in \GL(m,\mathbb{Z})}$ if $v \in
    \mathbb{Q}^m$.  In fact, write $(v_1,\ldots, v_m) = (a_1/d,
    \ldots, a_m/d)$.  Then the equation $(A-\id)v = l = (l_1,\ldots,
    l_m)$ is equivalent to $A (a_1,\cdots, a_m) = (a_1+l_1d,\ldots,
    a_m + l_kd)$. Now choose a prime $p$ not dividing any of the
    $a_i$ and $d$. We can solve the equations $a_i + l_id =0$ for $l_i$ in
    $\mathbb{Z}_p$ (since this is a field), which means that each $a_i
    + l_id$ will be divisible by $p$. But the equation $(A-\id)v = l$
    then implies that $(a_1,\ldots, a_m) = A^{-1}(a_1+l_1d,\ldots, a_m
    + l_kd)$ is divisible by $p$ as well.  So $l \notin \{Av-v, \; A
    \in \GL(m,\mathbb{Z})\}$.
\end{proof}

\begin{footnotesize}

\end{footnotesize}
\end{document}